\documentclass[12pt,a4paper]{article}
\usepackage{latexsym}
\usepackage{amssymb}
\expandafter\ifx\csname pre alinea.tex at\endcsname\relax \else
\endinput\fi \expandafter\chardef\csname pre alinea.tex
at\endcsname=\the\catcode`\@ \catcode`\@=11
\def\cases#1{\left\{\,\vcenter{\normalbaselines\m@th
    \ialign{$##\hfil$&\quad##\hfil\crcr#1\crcr}}\right.}
\def\matrix#1{\null\,\vcenter{\normalbaselines\m@th
    \ialign{\hfil$##$\hfil&&\quad\hfil$##$\hfil\crcr
      \mathstrut\crcr\noalign{\kern-\baselineskip}
      #1\crcr\mathstrut\crcr\noalign{\kern-\baselineskip}}}\,}
\def\pmatrix#1{\left(\matrix{#1}\right)}
\newdimen\p@renwd
\setbox0=\hbox{\rm B} \p@renwd=\wd0 
\def\bordermatrix#1{\begingroup \m@th
  \setbox\z@\vbox{\def\cr{\crcr\noalign{\kern2\p@\global\let\cr\endline}}%
    \ialign{$##$\hfil\kern2\p@\kern\p@renwd&\thinspace\hfil$##$\hfil
      &&\quad\hfil$##$\hfil\crcr
      \omit\strut\hfil\crcr\noalign{\kern-\baselineskip}%
      #1\crcr\omit\strut\cr}}%
  \setbox\tw@\vbox{\unvcopy\z@\global\setbox\@ne\lastbox}%
  \setbox\tw@\hbox{\unhbox\@ne\unskip\global\setbox\@ne\lastbox}%
  \setbox\tw@\hbox{$\kern\wd\@ne\kern-\p@renwd\left(\kern-\wd\@ne
    \global\setbox\@ne\vbox{\box\@ne\kern2\p@}%
    \vcenter{\kern-\ht\@ne\unvbox\z@\kern-\baselineskip}\,\right)$}%
  \null\;\vbox{\kern\ht\@ne\box\tw@}\endgroup}

\def\openup{\afterassignment\@penup\dimen@=}
\def\@penup{\advance\lineskip\dimen@
  \advance\baselineskip\dimen@
  \advance\lineskiplimit\dimen@}
\def\eqalign#1{\null\,\vcenter{\openup\jot\m@th
  \ialign{\strut\hfil$\displaystyle{##}$&$\displaystyle{{}##}$\hfil
      \crcr#1\crcr}}\,}
\newif\ifdt@p
\def\displ@y{\global\dt@ptrue\openup\jot\m@th
  \everycr{\noalign{\ifdt@p \global\dt@pfalse
      \vskip-\lineskiplimit \vskip\normallineskiplimit
      \else \penalty\interdisplaylinepenalty \fi}}}
\def\@lign{\tabskip\z@skip\everycr{}} 
\def\displaylines#1{\displ@y
  \halign{\hbox to\displaywidth{$\@lign\hfil\displaystyle##\hfil$}\crcr
    #1\crcr}}
\def\eqalignno#1{\displ@y \tabskip\centering
  \halign to\displaywidth{\hfil$\@lign\displaystyle{##}$\tabskip\z@skip
    &$\@lign\displaystyle{{}##}$\hfil\tabskip\centering
    &\llap{$\@lign##$}\tabskip\z@skip\crcr
    #1\crcr}}
\def\leqalignno#1{\displ@y \tabskip\centering
  \halign to\displaywidth{\hfil$\@lign\displaystyle{##}$\tabskip\z@skip
    &$\@lign\displaystyle{{}##}$\hfil\tabskip\centering
    &\kern-\displaywidth\rlap{$\@lign##$}\tabskip\displaywidth\crcr
    #1\crcr}}
\catcode`\@=\csname pre alinea.tex at\endcsname

\def\pa#1#2{\frac {\partial#1}{\partial #2}}
\def\parc#1{\frac {\partial}{\partial #1}}
\def\N{\mathbb N}
\def\R{\mathbb R}
\def\Z{\mathbb Z}
\def\C{\mathbb C}
\def\Q{\mathbb Q}

\def\Diff{\hbox{VF}}
\def\la{\lambda}
\def\comp{{\small\circ}}
\def\fle{\mathop{\longrightarrow}\limits}

\newtheorem{teo}{Theorem}[section]
\newtheorem{defi}{Definition}[section]
\newtheorem{lem}{Lemma}[section]
\newtheorem{coro}{Corollary}[section]
\newtheorem{prop}{Proposition}[section]

\def\proof{{\it Proof:\ }}
\def\fin{\hfill\hbox{$\Box$}\break}
\def\Cal#1{{\cal #1}}
\def\pmatrix#1{\left(\matrix{#1}\right)}

\def\dih#1{{\rm Diff}#1(\C^n,0)}
\def\dif#1{\widehat{\rm Diff}#1(\C^n,0)}
\def\camf#1{\hat\chi#1(\C^n,0)}
\def\dihdos#1{{\rm Diff}#1(\C^2,0)}
\def\difdos#1{\widehat{\rm Diff}#1(\C^2,0)}
\def\camfdos#1{\hat\chi#1(\C^2,0)}

\def\nor{\vartriangleleft}
\def\vect#1#2{\left({#1\atop #2}\right)}
\author{Fabio Enrique Brochero Mart\'{\i}nez\\ { \small Instituto de Matem\'atica Pura e Aplicada, IMPA}}
\date{
}
\title{\Large \bf Groups of germs of analytic diffeomorphisms in $(\C^2,0)$}

\begin{document}

\maketitle
\begin{abstract}
Let $\Cal G$ a group of germs of analytic diffeomorphisms in
$(\C^2,0)$. We find some remarkable properties  supposing that
$\Cal G$  is  finite, linearizable, abelian nilpotent and
solvable. In particular, if the group is abelian and has a generic
dicritic diffeomorphisms, then the group is a subgroup of a
1-parametric group. In addition, we study the topological
behaviour of the orbits of a dicritic diffeomorphisms. Last, we
find some invariants in order to know when two diffeomorphisms are
formally conjugates.
\end{abstract}

\section{Introduction}
The study of germs of holomorphic diffeomorphisms  and  finitely
generated groups of  analytic germs of diffeomorphisms fixing the
origin in one complex variable was started in the XIX century, and
 has been intensively studied by mathematicians
in the past century.
Such groups appear naturally when we study the holonomy group of
some leaf of codimension one holomorphic foliation.

From the Poincar\'e-Siegel linearization theorem, (See \cite{Ar} or
\cite{He}) it follows the study of local topology and the
analytical and topological classification of diffeomorphisms of
$(\C^n,0)$ having linear part in the Poincar\'e domain, or in the
Siegel domain that satisfy the Brjuno condition. The resonant case
in one dimension  is well known  too, the local topology and
topological classification was giving in Camacho \cite{Ca}.
Moreover the analytical classification in the resonant case  is
due to \'Ecalle \cite{Ec}, Voronin \cite{Vo}, Martinet and Ramis
\cite{MR} and Malgrange \cite{Ma}.

The topological classification  in dimension $ 2$ in the partially
hyperbolic case with resonances has been studied by Canille
\cite{Can} and the topological behavior in dimension $\ge 2$ in
the
 partial
hyperbolic case by Ueda \cite{Ue}.   In the case, tangent to the
identity,  Hakim \cite{Ha} and Abate \cite{Ab}, have shown the
existence of parabolic attractive points.

In addition, there  is an almost complete analytic classification
of the group of one dimensional germs when  we assume some
algebraic hypothesis as finiteness (See \cite{MM}), abelian,
solvable (see \cite{CM} and \cite{Na}). See \cite{EISV} for a
complete survey of this classification.

In this work, we deal with germs of diffeomorphisms in $(\C^2,0)$
and the finitely generated group of germs.

In  section 2 we give the  definitions and preliminary results.
%
In  section 3 we study the finite and linearizable  groups of
diffeomorphisms. We prove a generalization of Mattei-Moussu
topological criteria about finiteness of a group,

\

{\bf Theorem \ref{3.2}}{ \it Let $F\in \dihdos{}$. The group
generated by $F$ is finite if and only if  there exists a
neighborhood $V$ of 0, such that $|\Cal O_V(F,X)|<\infty$ for all
$X\in V$ and $F$ leaves  invariant infinite analytic varieties at
0.}

It is easy to see that  finite groups case is a particular case of
the linearizable case. Proposition \ref{p4.3} shows that in a
special case the topological linearization implies the analytical
linearization. In addition, for this special case we  construct
the moduli space (topological vs analytical) of the
diffeomorphisms that are conjugate that  diffeomorphism

{\bf Theorem \ref{4.1}} {\it  Let be $F\in \dih{}$, where
$A=DF(0)$ is a diagonalizable  matrix with norm 1 eigenvalue, i.e.
$A=diag(e^{2\pi i \lambda_1},\dots, e^{2\pi i \lambda_n})$ where
$\lambda_j\in \R$. Suppose that $F$ is topologically linearizable.
Then $$ \frac {\Cal H_{top}(F)}{\Cal H_{hol}(F)}\simeq {\pm}\frac
{SL(n,\Z)}{SL_A(n, \Z)}$$ where $SL_A(n,\Z)=\{B\in SL(n,\Z)|
(B-I)\lambda\in \Z^n\}.$ }

In particular,  we prove that the moduli space topological vs
$C^\infty$ is trivial.

Finally, we show that  a group of diffeomorphisms is linearizable
if and only if   there exists a  vector field with radial first
jet invariant by the  group-action, i.e.

\

{\bf Theorem \ref{4.2}} {\it A group $\Cal G\subset \dih{}$ is
analytically linearizable if and only if there exists a vector
field $\Cal X=\vec R+\cdots $, where $\vec R$ is a radial vector
field such $\Cal X$ is invariant for every  $F\in \Cal G$, i.e.
$F^*\Cal X=\Cal X$.}

\

In  section 4 we study the groups of diffeomorphisms supposing
that they have some algebraic structure. We prove that if $\Cal
G\subset \dihdos{}$ is a solvable group then its $7^{th}$
commutator subgroup is trivial. Furthermore, we characterize the
abelian subgroup of diffeomorphisms tangent to the identity, and
in the case when the group contains a dicritic diffeomorphism,
i.e. the group contains a diffeomorphism
$F(X)=X+F_{k+1}(X)+\cdots$ where $F_{k+1}(X)=f(X)X$ and $f$ is a
homogeneus polynomial of degree $k$,  we prove that the group is a
subgroup of a one parameter group. We write $\dihdos{_1}$ to
denote the group of diffeomorphisms tangent to the identity at
$0\in \C^2$.

\

{\bf Theorem \ref{5.5}} {\it  Let $\Cal G<\dihdos{_1}$ be  abelian
group, and $F\in \Cal G$ a dicritic diffeomorphism. Suppose that
$\exp(\mathfrak f)(x,y)=F(x,y)$ where $$\mathfrak
f=(f(x,y)x+p_{k+2}(x,y)+\cdots)\parc x +
(f(x,y)y+q_{k+2}(x,y)+\cdots)\parc y,$$ $f(x,y)$ is a homogeneous
polynomial of degree $k$ and
$g.c.d(f,xq_{k+2}(x,y)-yp_{k+2}(x,y))=1$, then $$\Cal G<\langle
\exp(t\mathfrak f)(x,y)|t\in \C\rangle.$$}

In  the section  5 we analyze the behavior of the orbits of a
diffeomorphism tangent to the  identity.  We prove a
generalization of the one dimensional flower theorem to two
dimensional  dicritic diffeomorphisms, i.e.

\

{\bf Theorem \ref{7.3}} {\it Let  $F:(\C^2,0)\to (\C^2,0)$ be a
dicritic diffeomorphism fixing zero,
 i.e. $F$ can be represented by a convergent series
$$F(x,y)=\left(\matrix{x+xp_k(x,y)+p_{k+2}(x,y)+\cdots\cr
 y+yp_k(x,y)+q_{k+2}(x,y)+\cdots}\right),$$
and $\tilde F=\Pi^*F:(\tilde\C^2,D)\to (\tilde\C^2,D)$ be the
continuous extension of the diffeomorphism  after  making the
blow-up in $(0,0)$. Then there exist open sets $U^+,U^-\subset
\tilde \C^2$ such that
\begin{enumerate}
\item $\overline{U^+\cup U^-}$ is a neighborhood of
$D\setminus\{(1:v)\in D|p_k(1,v)=0\}$. \item For all $P\in U^+$,
the sequence  $\{\tilde F^n(P)\}_{n\in \N}$ converge and
$\lim\limits_{n\to\infty}\tilde F^n(P)\in D$. \item For all $P\in
U^-$, the sequence  $\{\tilde F^{-n}(P)\}_{n\in \N}$ converge and
$\lim\limits_{n\to\infty}\tilde F^{-n}(P)\in D$.
\end{enumerate}}

\

An equivalent local theorem is proven  in the non dicritic case.

Finally, in  section 6 we show the formal classification of
diffeomorphisms tangent to the identity using the notion of the
semiformal conjugacy. We show that a representative diffeomorphism
found using semiformal conjugacy and  a cocycle determine its
formal conjugacy class.

In the case of dicritic diffeomorphisms, we find a rational
function that it is going to play  an equivalent role to the
residue in  a one dimension diffeomorphism.

\

{\bf Theorem \ref{8.3}} {\it Let $\tilde F\in \difdos{_1}$ be
dicritic diffeomorphism and
$F(x,v)=\left(\matrix{x+x^{k+1}p(v)+x^{k+2}(\cdots)\cr
v+x^{k+1}(\cdots)}\right)$ be the continuous extension of the
diffeomorphism  after  making the blow-up in $(0,0)$. Then there
exists a unique rational function $q(v)$ such that $F$ is
semiformally conjugate to
$$G_F=\left(\matrix{x+x^{k+1}p(v)+x^{2k+1}q(v)\cr v}\right)$$ in
$\overline \C\setminus \{p(v)=0\}$. In addition, $q(v)=\frac
{s(v)}{p(v)^{2k+1}}$ where $s(v)$ is a polynomial of degree
$2k+2+2k\partial(p(v))$.}

\

\section{Preliminaries}

Let  $\dif{}$ denote the  group of $n$-dimensional  formal
diffeomorphisms fixing zero, i.e.
$$\dif{}=\{H(X)=AX+P_2(X)+\cdots|A\in Gl(n,\C), P_i\in
\C^n[[X]]_i\},$$ where $X=(x_1,\dots,x_n)$ and $\C^n[[X]]_i$ is
the set of $n$-dimensional vectors with   coefficients homogeneous
polynomials of degree $i$.

Let  $\dih{}\subset \dif{}$ denote the pseudo-group of germs of
holomorphic diffeomorphisms, i.e. the power series that represent
every   $H\in \dih{}$ convergent in some neighborhood of $0\in
(\C^n,0)$.

For each $j\ge 2$, let  $\dif{_j}$ (resp. $\dih{_j}$) denote the
subgroup of formal (resp. analytic)  diffeomorphisms $j$-flat,
i.e. $F(X)\in\dif{_j}$ (resp. $F(X)\in\dih{_j}$) then
$F(X)=X+P_j(X)+\cdots $.

The same way define $\camf{_j}$ the Lie algebra of formal vector
fields of $\C^n$, $j$-flat in $0\in \C^n$, i.e.
$$\camf{_j}=\{F_1(X)\parc {x_1}+\cdots+ F_n(X)\parc {x_n}|\hbox{
where } F_k\in \mathop{\oplus}_{i=j}^\infty\C[[X]]_i\}.$$

\begin{prop}\label{expo} The exponential map $\exp:\camf{_j}\to \dif{_j}$ is a bijection
for every $j\ge 2$.
\end{prop}\label{exp}

\proof The proof follows from a  straightforward generalization of
the proposition 1.1.  in \cite{MR}.\fin

Notice  that $\camf{_j}$ is the formal Lie algebra associate to
the formal Lie group $\dif{_j}$. It is not difficult to prove that
if the vector  field is holomorphic, then  the associated
diffeomorphism
 is also holomorphic. The converse is false in general.

\begin{defi} $\Cal H, \Cal G\subset \dif{}$ are called
{\bf formally (analytically) conjugate} if there exists $g\in
\dif{}$ ($g\in \dih{}$) such that $g\comp \Cal H\comp
 g^{-1}=\Cal G$. \end{defi}

\begin{defi} $\Cal H, \Cal G\subset \dih{}$ are called
{\bf topologically  conjugate} if there exists $t:\Cal H\to \Cal
G$ bijective group homomorphism and   $g:U\to g(U)$ local
homeomorphism at 0, such that $F\circ g=g\circ t(F)$ for all $F\in
\Cal H$ in some neighborhood of 0. \end{defi}

In particular, some group $\Cal G$ is called linearizable if there
exists $g$ such that $g\circ \Cal G\circ g^{-1}$ is a linear group
action.   Moreover, if $g$ is a diffeomorphism formal, analytic or
continuous then the group is called formally, analytically or
topologically linearizable respectively. Let $\Lambda_{\Cal G}$
denote the set $\{D G(0)| G\in \Cal G\}$. To understand when a
diffeomorphism is formally and analytically linearizable we need
the following definition

\begin{defi}Let $G\in \dih{}$. $G$ is called {\bf resonant} if the
 eigenvalues of $G'(0)$,   $\la:=(\la_1,\dots,\la_n)$, satisfies
some relation like $$\la^m-\la_j=\la_1^
{m_1}\cdots\la_n^{m_n}-\la_j=0$$
 for some
$j=1,\dots,n$ and $m\in \N^n$ where $|m|=m_1+\cdots+m_n\ge2$.
Otherwise,  it is called non resonant.
\end{defi}

\begin{defi} Let $A\in Gl(n,\C)$ and suppose that $A$ is not
resonant. We say that $A$ satisfies the Brjuno condition if $$
\sum_{k=0}^\infty 2^{-k} \log (\Omega^{-1} (2^{k+1}))<+\infty$$
where $\Omega(k)=\inf\limits_{{2\le |j|\le k \atop 1\le i\le n}} |
\lambda^j-\lambda_i|$.
\end{defi}

\begin{teo} Let $G\in \dih{}$, where $G(X)=AX+\cdots$, A
diagonalizable, non resonant. Then $G$ is formally linearizable.
In addition, if $G$ satisfies the Brjuno condition,  then this
linearization is in fact holomorphic.
\end{teo}

In the resonant case in dimension one  is easy to see that
$f(x)=g\circ h(x)$ where $g$ and $h$ commute, $g^n=id$ and $h$ is
tangent to the identity. Notice that $g$ and $h$ are in general,
formal diffeomorphisms.

\begin{teo}[Camacho \cite{Ca}] Let $h(x)=x+a x^{k+1}+\cdots$, then $h$ is topologically
conjugate with $x\mapsto x+x^{k+1}$.
\end{teo}
Thus, the topological classification is very simple. In the same
way the formal classification is simple
\begin{teo}
Let $h(x)= x+ax^{k+1}+ \cdots$, then there exists $\rho$ such that
$h$ is formally conjugate with $f_{a,k,\rho}(z)=\exp(a\frac
{z^{k+1}}{1+\rho z^k}\parc z)$.
\end{teo}
$\rho(h)$ is called the residue of $h$ and
$\rho(h)=\displaystyle{\frac 1{2\pi i}\oint \frac 1{h(x)-x}dx}$.
In particular, $h$ is formally conjugate with $x\mapsto
x+x^{k+1}+\rho x^{2k+1}$.

In addition, there exists a unique $g_h$ tangent to the identity
such that $g_h^*h=g_h\circ h\circ g_h^{-1}=\exp(a\chi_{k,\rho})$
where $\chi_{k,\rho}= \frac {z^{k+1}}{1+\rho z^k}\parc z$. This
$g_h$ is called the {\bf normalizer transformation} of $h$.

\section{Finite Groups}

The  finite group of germs of holomorphic diffeomorphisms in one
dimension appears naturally in the study of the holonomy group  of
the  local foliation with holomorphic first integral around a
singular point.

\begin{prop} Let  $\Cal H$ be a finite subgroup of\ $\dif{}$
(resp. $\dih{}$) then $\Cal H$ is formally (analytic)
linearizable, and it is isomorphic to a finite subgroup of
$Gl(n,\C)$. \end{prop}

\proof Let $g(X)=\sum\limits_{H\in \Cal H} (H'(0))^{-1}H(X)$. Note
that $g$ is a diffeomorphism because $g'(0)=\#(\Cal H) I$, and
moreover for all $F\in \Cal H$ $$\eqalign{g\comp F(X)&=\sum_{H\in
\Cal H} (H'(0))^{-1}H\comp F(X) 
=F'(0)\sum_{H\in \Cal H} ((H\comp F)'(0))^{-1}H\comp F(X)\cr
&=F'(0)g(X) }$$ Thus $g\comp F\comp g^{-1}(X)=F'(0)X$.  In fact,
we obtain a injective groups homomorphism $$\eqalign{\Cal
H&\fle^{\Lambda} Gl(n,\C)\cr F&\fle (g\comp F\comp
g^{-1})'(0).}\eqno{\ \atop\fin}$$

Denote by $\Lambda_{\Cal H}\subset Gl(n,\C)$ the group of  linear
parts of the diffeomorphisms in $\Cal H$,  and $p=\#(\Cal H)$,
then for all matrix $A\in \Lambda_{\Cal H}$, $A^p=I$. We claim
that $A$ is diagonalizable, in fact, for the  Jordan canonical
form  theorem there exists $B\in Gl(n,\C)$ such that
$BAB^{-1}=D+N$ where  $D$ is a   diagonal matrix and $N$ is such
that $ND=DN$ and $N^l=0$ for some $l\in \N$,  but $A^p=I$  then
$(D+N)^p=I$. To prove the claim we need the following lemma

\begin{lem} Let  $l$ be the least integer such that $N^l=0$, then $I,N,\dots,
N^{l-1}$ are $\C$-linearly independent. \end{lem} \proof Let
$a_mN^m+\cdots+a_{l-1}N^{l-1}=0$  be a linear combination  with
$a_m\ne 0$  $m\ge 0$,  multiply by $N^{l-m-1}$ we obtain that
$a_mN^{l-1}=0$, and for the  minimality of $l$ we conclude that
$a_m=0$, it is a contradiction. \fin

\quad We can suppose, that   each block of the Jordan canonical
form is of the form   $D=\la I$, thus $(D+N)^p=(\la I+N)^p=I$,
therefore $(\la^p-1)I+p\la^{p-1}N+\cdots +N^p= 0$,  for the lemma
we have that $\la^p=1$ and $N=0$.

In order to show a topological criteria to know when a
diffeomorphism is finite we need the following modification of the
Lewowicz lemma.

\begin{lem} Let $M$, $0\in M$, be a complex analytic variety of $\C^n$ and
$K$ a connected component of 0 in $\overline B_r(0)\cap M$.
Suppose that $f$ is a homeomorphism from $K$ to $f(K)\subset M$
such that $f(0)=0$. Then there exists $x\in \partial K$ such that
the number of iterations  $f^m(x)\in K$ is infinity.
\end{lem}
\proof Denote by $\overline \mu=\mu|_K$ and
$\mu=\mu|_{\mathop{K}\limits^\circ}$ the number of iteration in
$K$ and $\mathop{K}\limits^\circ$. It is easy to see that
$\overline\mu$ is upper semicontinuous, $\mu$ is under
semicontinuous and $\overline\mu(x)\ge  \mu(x)$ for all $x\in
\mathop{K}\limits^\circ$. Suppose by contradiction that $\overline
\mu(x)<\infty$  for all $x\in \partial K$, therefore exists $n\in
\N$ such that $\overline \mu(x)<n$ for all $x\in
\partial K$.
Let  $A=\{x\in K|\overline \mu(x)<n\}\supset \partial K$ and
$B=\{x\in \mathop{K}\limits^\circ| \mu(x)\ge n\}\ni 0$ open set,
and $A\cap B=\emptyset$ since $\overline\mu(x)\ge \mu(x)$.

Using the fact that $K$ is a connected set, there exists $x_0\in
K\setminus(A\cup B)$ i.e $\overline\mu(x_0)\ge n>\mu(x_0)$, then
the orbit of $x_0$ intersects the border of $K$, which is a
contradiction since $\partial K\subset A$ implies $x_0\in A$.\fin

Let $f:U\to f(U)$ be a homeomorphism with $f(0)=0$ and $x\in U$.
We denote by $\Cal O_U(f,x)$ the $f$-orbit of $x$ that do not
leave $U$, i.e. $y\in \Cal O_U(f,x)$ if and only if
$\{x,f(x),\dots, y=f^{[k]}(x)\}\subset U$ or
$\{x,f^{[-1]}(x),\dots, y=f^{[-k]}(x)\}\subset U$ for some $k\in
\N$.

\begin{teo}\label{3.2} Let $F\in \dihdos{}$. The group generated by $F$ is finite if and
only if  there exists a neighborhood $V$ of 0, such that $|\Cal
O_V(F,X)|<\infty$ for all $X\in V$ and $F$ leaves  invariant
infinite analytic varieties at 0.

\end{teo}

\proof $(\Rightarrow)$ Let $N=\#\langle F\rangle $ and $h\in
\dihdos{}$ such that $h\comp F\comp h^{-1}(x,y)=(\lambda_1
x,\lambda_2 y)$ where $\lambda_1^N=\lambda_2^N=1$. It is clear
than $|\Cal O(F,X)|\le N$ for all $X$ in the domain of $F$, and
$M_c=\{h(x,y)|x^N-cy^N=0\}$ is a complex analytic variety
invariant by $F$ for all $c\in\C$.

$(\Leftarrow)$ Without loss of generality we suppose that
$V=\overline B_r(0)$ where $F(V)$ and $F^{-1}(V)$ are well
defined. Let $M$ be a  $F$-invariant complex analytic variety and
$K_M$ the connected component of 0 in $M\cap V$. Let $A_1=K$,
$A_{j+1}=K\cap F^{-1}(A_j)$ and $C_n$ the connected component of 0
in $A_n$. It is clear, by construction that $A_n$ is the set of
point of $K_M$ with $n$ or more iterates in $K_M$. Moreover, since
$A_n$ is compact and $C_n$ is compact and connected, it follows
that $C_M=\cap_{j=1}^\infty C_n$ is compact and connected too, and
therefore $C_M=\{0\}$ or $C_M$ is non enumerable.

We claim that $C_M\cap \partial K\ne \emptyset$ and then this is
non enumerate. In fact, if $C_M\cap \partial K= \emptyset$ there
exists $j$ such that $C_j\cap \partial K=\emptyset$. Let $B$
compact connected neighborhood of $C_j$ such that $(A_j\setminus
C_j)\cap B=\emptyset$, therefore for all $X\in \partial B$ we have
$\Cal O_K(F,X)<N$, that  is a contradiction by the lemma.

In particular, $C=C_{\C^2}$ is a set of  point with  infinite
orbits in $V$ and therefore every point in $C $ is  periodic. If
we denote $D_n=\{X\in C| F^{n!} (X)=X\}$, it is clear that $D_n$
is a close set and $D_n\subset D_{n+1}$, moreover
$C=\cup_{n=1}^\infty D_n$, then exists $n\in \N$ such that
$C=D_n$. Let $G=F^{n!}$ where it is  well defined, observe  that
$C$ is in the domain U of $G$ and $C\subset\{X\in U|G(X)=X\}=L$.
Since $L$ is a complex analytic variety of $U$ that contain $C$
then its dimension is 1 or 2. The case $\dim L=1$ is impossible
because $C_M\subset C\subset L$ for all $M$ analytic variety
$F$-invariant, contradicting  that fact that $\Cal O_2$ is
Noetherian ring. In the case $\dim L=2$ follows that $F^{n!}(X)=X$
for all $X\in U$, therefore $\langle F\rangle$ is finite.\fin

This theorem can be extended to $\C^n$ if we suppose that there
exist infinitely many one dimensional invariant analytic varieties
at 0 in general position.

\begin{defi}  Let $\Cal G<\dih{}$ be a group of germs of diffeomorphisms.
\begin{itemize}
\item[a)] $\Cal G$
 is called {\bf periodic} if every element $F\in \Cal G$
has finite order.
\item[b)] $\Cal G$ is called {\bf locally finite} if every subgroup finitely
generated is finite.
\item[c)] $\Cal G$ is called {\bf a group of exponent $d$}, if $F^{[d]}=id$, for every element
$F\in \Cal G$.
\end{itemize}
\end{defi}

\begin{prop}\label{p3.2} Let $\Cal G$ be a periodic group, then  the
homomorphism \break $\eqalign{\Cal G&\fle^\Lambda Gl(n,\C)\cr
G&\longmapsto G'(0),}$ is injective.
\end{prop}
\proof In fact, suppose  by  contradiction that there exists $F\in
\Cal G$ such that $\Lambda(F(X))=I$, i.e., $F(X)=X+P_k(X)+\cdots$,
then, by straightforward calculation  we have that
$F^{[r]}(X)=X+rP_k(X)+\cdots\ne Id$,  therefore, the unique
element of $\Cal G$ tangent to the identity is  itself. \fin

\begin{coro} Let  $\Cal G<\dih{}$ such that, every element $F\in \Cal G$ has
finite order, then $\Cal G_1\le \Cal G$ the subgroup of element
tangent to the identity is trivial, i.e., $\Cal G_1=\{id\}$.
\end{coro}

\begin{prop} Let $\Cal G<\dih{}$
\begin{itemize}
\item[a)] If there exists $p\in \N^*$ such that, for every $F\in \Cal G$,
$F^p=id$, then $\Cal G$ is finite.
\item[b)] If $\cal G$ is periodic  then every finitely generated subgroup is
finite.
\end{itemize}
\end{prop}
For the proposition above we only need to prove that the item a)
and b) are true  for the group $\Lambda_{\Cal G}<Gl(n,\C)$. This
fact follows from the next theorems
\begin{teo}[Burnside] 
If $G<Gl(n,\C)$ is a group with finite exponent $m$, then
$|G|<m^{n^3}$.
\end{teo}

\begin{teo}[Schur]
If $G<Gl(n,\C)$ is a periodic group then $G$ is locally finite.
\end{teo}

\begin{coro}Let $\Cal G\in \dihdos{}$ be a finite generate group. The group  $\Cal G$ is finite if and
only if  there exists a neighborhood $V$ of 0, such that $|\Cal
O_V(\Cal G,X)|<\infty$ for all $X\in V$ and $F$ leaves  invariant
infinite analytic varieties at 0.
\end{coro}

\proof From the theorem \ref{3.2} we know that every element of
$\Cal G$ has finite order. Therefore from the proposition
\ref{p3.2} follows that $\Cal G$ is isomorphic to the group $\{D
G(0)|G\in \Cal G\}$. Finally from the Schur and Burnside theorem
we conclude that $\Cal G$ is finite.

\begin{prop}\label{p4.3} Let $F\in \dih{}$ such that $A=DF(0)$ is a diagonalizable matrix
 and the
eigenvalues have norm 1. Then $F$ is analytically linearizable if
and only if $F$ is topologically linearizable.
\end{prop}
\proof Suppose that $F$ is topologically linearizable and $A$ is a
diagonal matrix,  i.e. there exists a homeomorphism $h:U\to h(U)$
where $U$ is a neighborhood of $0\in \C^n$, such that $h\circ
F(X)=Ah(X)$ where $X\in F^{-1}(U)=W$.

Claim: there exists a neighborhood $V\subset W$ of $0\in \C^n$
such that $F(V)=V$. In fact, let $r>0$ such that $B(0,r)\subset
h(U)\cap W$ and $V=h^{-1}(B(0,r))$, it is clear that
$$F(V)=h^{-1}(Ah(V))=h^{-1}(A(B(0,r)))=h^{-1}(B(0,r))=V.$$

Since $H_m(X)=\frac 1m\sum\limits_{j=0}^{m-1} A^{-j} F^j(X)$ is a
normal family of holomorphic diffeomorphisms defined from   $V$ to
the smaller Reinhardt domain that contains $V$, then
$$H(X)=\lim_{m\to \infty} \frac 1m\sum\limits_{j=0}^{m-1} A^{-j}
F^j(X)$$ is  holomorphic in $V$, with $H'(0)=Id$, it follows that
$H$ is a local diffeomorphism at $0\in \C^n$.

Moreover, notice that $H_m(F(X))=AH_m(X)+\frac 1m
(A^{m-1}F^m(X)-X)$ and $$\lim\limits_{m\to \infty}\frac
{|A^{m-1}F^m(X)-X|}m=0.$$ We conclude therefore that
$H(F(X))=AH(X)$ for all $X\in V$. \fin

\begin{defi} For all $F\in \dih{}$,  let $$\Cal H_{top}(F)=\{G\in \dih{}|
G\hbox{ is topologically conjugate with } F\}$$ and $$\Cal
H_{hol}(F)=\{G\in \dih{}| G\hbox{ is holomorphically  conjugate
with } F\}.$$ \end{defi}

\begin{teo}\label{4.1} Let be $F\in \dih{}$, where $A=DF(0)$ is a diagonalizable  matrix
with norm 1  eigenvalue, i.e. $A=diag(e^{2\pi i \lambda_1},\dots,
e^{2\pi i \lambda_n})$ where $\lambda_j\in \R$. Suppose that $F$
is topologically linearizable.  Then $$ \frac {\Cal
H_{top}(F)}{\Cal H_{hol}(F)}\simeq {\pm}\frac {SL(n,\Z)}{SL_A(n,
\Z)}$$ where $$SL_A(n,\Z)=\{B\in SL(n,\Z)| (B-I)\lambda\in
\Z^n\}.$$
\end{teo}
\proof

From the proposition \ref{p4.3}, we know that every element of
$\Cal H_{top}(F)$ is analytically linearizable. Let $G\in \Cal
H_{top}(F)$, $H$ local homeomorphism such that $F\circ H=H\circ
G$, $h:=H|_{S^1\times\cdots{\times}S^1}$,
$f=F|_{S^1\times\cdots{\times}S^1}$ and $g=G|_{S^1\times
\cdots{\times}S^1}$. If $h=(h_1,\dots,h_n)$,  since the norm of
$h_j$  is invariant by $F $ and $G$, we can suppose without loss
of generality that $h(S^1\times\cdots\times
S^1)=S^1\times\cdots{\times}S^1$ and then the following diagram
commute $$\matrix{ S^1\times\cdots\times
S^1&\fle^h&S^1\times\cdots{\times}S^1\cr \ \downarrow^f&&\
\downarrow^g\cr S^1\times\cdots\times
S^1&\fle^h&S^1\times\cdots{\times}S^1. } $$ Let $\Pi:\R^n\to
S^1\times\cdots{\times}S^1$, \quad $(z_1,\dots,z_n)\mapsto (e^{2\pi
iz_1},\dots, e^{2\pi i z_n})$ be  a universal covering of
$S^1\times\cdots{\times}S^1$,

Let  $C$  be a lifting of $H|_{S^1\times\cdots{\times}S^1}$.
Observe that $C|_{\Z^n}:\Z^n\to \Z^n$ can be represented as a
linear transformation, and since the volume of $C([0,1]^n)$ is 1,
then $C|_{\Z^n}=M\in{\pm}SL(n,\Z)$. Therefore, $C(Z)=M([\![
Z]\!])+\theta(\{Z\})$, where $[\![ Z]\!]=([\![ z_1]\!],\dots,[\![
z_n]\!])$, $\{Z\}=Z-[\![ Z]\!]$ and $\theta=C|_{[0,1)^n}$. In
addition, $$C(Z)+\mu+C(Z+\lambda)\in \Z^n$$ where
$\mu=(\mu_1,\dots, \mu_n)$ and $DG(0)=diag(e^{2\pi i \mu_1},\dots,
e^{2\pi i \mu_n})$. In particular, making  $Z=0$ we have
$$C(j\lambda)-j\mu\in \Z^n,\quad\hbox{for all $j\in \Z$}$$ and
replacing $\mu$ by $\mu+m$ for some $m\in \Z^n$, we can suppose
that $C(j\lambda)=j\mu$. It follows that
$$M([\![j\lambda]\!])-[\![j\mu]\!]=\{j\mu\}-\theta(\{j\lambda\}).$$
From the continuity of $\theta$ and using that there exists a
positive integers increasing sequence  $\{a_n\}_{n\in \N}$ such
that $\{a_n\lambda\}\mathop{\longrightarrow}\limits_{n\to \infty}
0$ and $\{a_n\mu\}\mathop{\longrightarrow}\limits_{n\to \infty}
0$, we have that $M([\![a_n\lambda]\!])=[\![a_n\mu]\!]$ for all
$n\gg 0$, therefore $$M\lambda=\lim_{n\to\infty} \frac 1{a_n}
M([\![a_n\lambda]\!])=\lim_{n\to\infty} \frac
1{a_n}[\![a_n\mu]\!]=\mu.$$

Thus, if we have two homeomorphism that conjuge $F$ and $G$, each
one has associated a matrix $M_1,M_2\in SL(n,\Z)$ such that
$M_1\lambda-\mu\in \Z^n$ and $M_2\lambda-\mu\in \Z^n$, then
$(M_1M_2^{-1}-I)\lambda\in \Z^n$, as we want to prove. \fin

 Finally, for each element $M\in{\pm}\frac
{SL(n,\Z)}{SL_A(n,\Z)}$ it is easy to see that a linear
holomorphic diffeomorphism $G$ $$(x_1,\dots,x_n)\mapsto (e^{2\pi i
\mu_1}x_1,\dots,  e^{2\pi i \mu_n}x_n)\quad\hbox{where}\quad
\mu\equiv M \lambda \pmod 1,$$ and the $C^\infty$- diffeomorphism
$H$
$$H(x,y)=(|x_1|e^{m_{11}\arg(x_1)+\cdots+m_{1n}\arg(x_n)},|x_n|e^{m_{n1}\arg(x_1)+\cdots+m_{nn}\arg(x_n)})$$
satisfies $H\circ F=G\circ H$.\fin

\begin{coro} Let $F,G\in \dih{}$ such that they are topologically conjugates,
 where $A=DF(0)$ is a diagonalizable  matrix with norm 1
eigenvalue, then they are $C^{\infty}$ conjugate.
\end{coro}

\proof Let $h_1, h_2$ holomorphism diffeomorphisms such that
$F_1=h_1\circ F\circ h_1^{-1}$ and $G_1=h_2\circ G\circ h_2^{-1}$
are linear transformations. From the theorem there exists a matrix
$M=(m_{ij})\in {\pm} SL(n,\Z)$ such that the
$C^{\infty}$-diffeomorphism $H(X)=\Bigl(|x_1|e^{\sum_{j=1}^n
m_{1j}\arg(x_j)},\dots,|x_n|e^{\sum_{j=1}^n
m_{nj}\arg(x_j)}\Bigr)$ conjugate $F_1$ and $G_1$. Therefore
$h_2^{-1}\circ H\circ h_1$ conjugate $F$ and $G$.\fin

Now observe that the radical vector field is invariant by the
action of every linear diffeomorphism. This fact characterizes
every linearizable group

\begin{teo}\label{4.2} $\Cal G\subset \dih{}$ is a group analytically linearizable if and only if
there exists a vector field $\Cal X=\vec R+\cdots $, where $\vec
R$ is a radial vector field such $\Cal X$ is invariant for every
$F\in \Cal G$, i.e. $F^*\Cal X=\Cal X$.
\end{teo}

\proof $(\Rightarrow)$ Suppose that $\Cal G$ is linearizable, i.e.
there exists $g:(\C^n,0)\to (\C^n,0)$ such that $g^{-1}\circ \Cal
G\circ g=\{DF(0)Y|F\in \Cal G\}$. Since $(AY)^* \vec R=\vec R$ for
all $A\in Gl(2,\C)$, in particular for every element $F\in \Cal G$
we have $$\eqalign{ \vec  R&=(g^{-1}\circ F\circ g)^*\vec R=
D(g^{-1}\circ F\circ g)_{(g^{-1}\circ F^{-1}\circ g(Y))}{\cdot}\vec
R (g^{-1}\circ F^{-1}\circ g(Y))\cr &=Dg^{-1}_{(g(Y))}\cdot
DF_{(F^{-1}(g(Y))}{\cdot}Dg_{((g^{-1}\circ F^{-1}\circ g(Y))}\cdot
(g^{-1}\circ F^{-1}\circ g(Y)) } $$ Replacing $X=g(Y)$ and
multiplying by $D g_{g^{-1}(Y)}$ we have  that
$$DF_{(F^{-1}(X)}{\cdot}Dg_{((g^{-1}\circ F^{-1}(X))}\cdot
(g^{-1}\circ F^{-1}(X))=Dg_{g^{-1}(Y)} g^{-1}(X),$$ i.e. denoting
$\Cal X=Dg_{g^{-1}(X)} g^{-1}(X)$ we have that $F^*\Cal X=\Cal X$.
It is easy to see that $\Cal X=\vec R+\cdots$.

$(\Leftarrow)$ Now suppose that  $F^*\Cal X=\Cal X$ for every
element $F\in \Cal G$, where $\Cal X=\vec R+\cdots$.

Since every eigenvalue of the linear part of $\Cal X$ is 1, then
$\Cal X$ is in the Poincar\'e domain without resonances, therefore
there exists a analytic diffeomorphism $g:(\C^n,0)\to (\C^n,0)$
such that $g^*\Cal X=\vec R$, i.e. $\Cal X=(Dg_{(X)})^{-1}\cdot
g(X)$.

We claim that $g\circ F\circ g^{-1}(Y)=DF(0)Y$  for every $F\in
\Cal G$. In fact, from the same procedure as before we can observe
that $$(g\circ F\circ g^{-1})^*\vec R=\vec R.$$ Now, if we suppose
that $g\circ F\circ g^{-1}=AX+P_l(X)+P_{l+1}(X)+\cdots$, where
$P_j(X)$ is a polynomial vector field of degree $j$, then it is
easy to see using the Euler equality that $$(g\circ F\circ
g^{-1})^*\vec R=AX+lP_l(X)+(l+1)P_{l+1}(X)+\cdots,$$ and therefore
$P_j(X)\equiv 0$ for every $j\ge 2$. \fin

\section{Groups with some algebraic structure}

In this section we are going to study the groups of
diffeomorphisms, when they have some additional algebraic
structure as abelian, nilpotent or  solvable.

Let $\Cal G<\dih{}$. The upper central series $$Z_0=\{id\}\subset
Z_1(\Cal G)\subset\cdots \subset Z_n(\Cal G)\subset\cdots$$ of
$\Cal G$ is defined inductively where $Z_{j+1}(\Cal G)/Z_j(\Cal
G)$ is the center of $\Cal G/Z_j(\Cal G)$. The group $Z_j(\Cal G)$
is called the  {\bf $j$-th hypercenter} of $\Cal G$. $\Cal G$ is
called {\bf nilpotent} if $\Cal G=Z_l (\Cal G)$ for some $l$.  The
smallest $l$ for which $\Cal G=Z_l(\Cal G)$ is the {\it nilpotency
class} of $\Cal G$.

In the same way  for any  subalgebra $\Cal L\subset \camf{}$, we
define $Z_j[\Cal L]$ the $j$th hypercenter of $\Cal L$ inductively
such that $Z_{j+1}(\Cal L)/Z_j(\Cal L)$ is the center of $\Cal
L/Z_j(\Cal L)$ and we say that $\Cal L$ is nilpotent of class $l$
if $l$ is the smallest integer  such that $Z_l(\Cal L)=\Cal L$.

The commutator series $$\Cal G^0=\Cal G\supset\Cal
G^1\supset\cdots\supset\Cal G^j\supset\cdots$$ is defined
inductively where  $\Cal G^j=[\Cal G^{j-1},\Cal G^{j-1}]$ is the
{\bf $j$-th commutator subgroup}. $\Cal G$ is called  {\bf
solvable}, if there exists a  positive integer $l$ such that
$\Cal G^l=\{Id\}$. It is obvious that every nilpotent group is
solvable.

Denoting $\Cal G_1\nor \Cal G$ the normal subgroup of  the
diffeomorphisms tangent to the identity, it is easy to see that
$$\Cal G/\Cal G_1\sim\Lambda_{\Cal G}=\{DG(0)| G\in \Cal G\},$$
therefore $\Cal G$ is solvable if and only if $\Lambda_{\Cal G}$
is solvable and $\Cal G_1$ is solvable.

In addition, if   $\Cal G$ is solvable, and let $\{id\}=\Cal
G^l\nor  \Cal G^{l-1}\nor \cdots\nor \Cal G^1 \nor \Cal G$ the
resolution string, then using the group homomorphism $\Lambda$ we
obtain a new resolution string $$\matrix{\Cal G^l\nor & \Cal
G^{l-1}\nor &\cdots& \nor &\Cal G^1 \nor &\Cal G\cr
\downarrow^\Lambda&\downarrow^\Lambda&\vdots&&\downarrow^\Lambda
&\downarrow^\Lambda\!\!\cr G^l\nor &G^{l-1}\nor &\cdots& \nor &G^1
\nor &G}$$ where $G^j=\Lambda(\Cal G^j)$. Denote by $height(\Cal
G)=l$ the height of the resolution string of $\Cal G$,  then it is
clear   $height(G)\le height(\Cal G)$. Note that in the case where
$G< Gl(n,\C)$ is a linear solvable group, it is known that
$height(G)$ is limited by a function that only  depends  on $n$,
$\rho(n)$. In fact  Zassenhaus proved that $\rho(n)<2n$ and Newman
(see \cite{Ne}) has found $\rho(n)$ for all $n$, in particular,
$\rho(2)=4$ and $\rho(3)=5$.

Let $\Cal G_j=\{G\in \Cal G_1| \hbox{where $G$ is $k$-flat}, k>
j\}$. It is easy to see that   $$\cdots\nor\Cal G_3\nor \Cal
G_2\nor \Cal G_1$$ is a normal series and $\Cal G_j/\Cal G_{j+1}$
are  abelian groups.

If $\Cal G_1 \cap \dif{_k}=\{id\}$ for some $k\in \N$ the
following lemma show that $\Cal G_1$ is solvable, and provide a
necessary condition in order to two diffeomorphisms tangent to the
identity commute.

\begin{lem}\label{comuta} Let $F\in \dif{_{r+1}}$ and $G\in \dif{_{s+1}}$,
then $$F(G(X))-G(F(X))\!=\nabla F_{r+1}(X){\cdot}G_{s+1}(X) -\nabla
G_{s+1}(X){\cdot}F_{r+1}(X)+O(|X|^{r+s+2}).$$ In particular,
$[F,G]\in \dif{_{r+s+1}}$
\end{lem}
\proof Let
$F(X)=X+\sum\limits_{k=r}^{r+s}F_{k+1}(X)+O(|X|^{r+s+2})$ and
$G(X)=X+\sum\limits_{k=s}^{r+s}G_{j+1}(X)+O(|X|^{r+s+2})$, then

\begin{equation}
\eqalign{F(G(X))&=X+\sum_{k=s}^{r+s}G_{j+1}(X)+O(|X|^{r+s+2})+\cr
&\qquad\quad+\sum_{k=r}^{r+s}F_{k+1}\Bigl(
X+\sum_{k=s}^{r+s}G_{j+1}(X)
+O(|X|^{r+s+2})\Bigr)+O(|X|^{r+s+2})\cr
&=X+\sum_{k=s}^{r+s}G_{j+1}(X)+
\sum_{k=r}^{r+s}\Bigl(F_{k+1}(X)+\nabla F_{k+1}(X)\cdot
G_{s+1}(X)\cr &\hskip 6.3cm
+O(|X|^{k+s+2})\Bigr)+O(|X|^{r+s+2})\cr
&=X+\sum_{k=s}^{r+s}G_{j+1}(X)+ \sum_{k=r}^{r+s}F_{k+1}(X)+\nabla
F_{r+1}(X){\cdot}G_{s+1}(X)+O(|X|^{r+s+2})}
\label{comuta1}\end{equation} In the same way we have
\begin{equation}
G(F(X))=X+\sum_{k=s}^{r+s}G_{j+1}(X)+
\sum_{k=r}^{r+s}F_{k+1}(X)+\nabla G_{s+1}(X)\cdot
F_{r+1}(X)+O(|X|^{r+s+2}). \label{comuta2}
\end{equation}
The lemma follows subtracting the equations (\ref{comuta1}) and
(\ref{comuta2}).\fin

\begin{coro} Let $\Cal G$ be a group of diffeomorphisms and suppose that
there exists $l\in \N$ such that $\Cal G\cap \dih{_l}=\{id\}$.
Then $\Cal G$ is solvable if and only if $\Cal G^{\rho(n)}\subset
\dih{_1}$, where $\rho(n$) is the Newman  function.
\end{coro}

\begin{teo}
Let $\Cal G$ be a solvable subgroup of $\dif{}$, then $\Cal
G^{\rho(n)+1}$ is a nilpotent group.
\end{teo}

\proof Observe that If $\Cal G$ is a solvable group, then $\Cal
H=\Cal G^{\rho (n)}\subset \dih{_1}$ is a solvable group. Let
$\Cal L$ be the algebra associate to the group $H$ by  the $\exp$
function. For every integer $i\ge 0$ denote $\Diff_i=\frac
{\camf{}}{\camf{_i}}$. Notice that for all $j>i$ there exists a
natural projection $\Diff_j\to \Diff_i$, therefore we can think of
$\camf{}$ as a projective limit $\lim\limits_{j\to \infty}
\Diff_j$, moreover $\Diff_j$ has a natural structure of finite
dimensional complex algebraic  linear algebra. Let $[\Cal L]_i$ be
the projection of $\Cal  L$ in $\Diff_i$,  and $\overline{[\Cal
L]_i}$ the algebra of every linear combination . It follows
 that
  $\overline{[\Cal L]_i}$ is  a
connected solvable  Lie algebra  of the same height as $\Cal L$.
Since $\overline{[\Cal L]_i}$ is isomorphic to a linear Lie
algebra, we know from the Lie-Kolchi theorem that they can be
represented by upper triangular matrices, then $[\overline{[\Cal
L]_i},\overline{[\Cal L]_i}]$ is represented by nilpotent
triangular matrices for every $i\in \N$. It follows that $\Cal
L^1=\lim\limits_{j\to \infty} [\Cal L^1]_j$ is a nilpotent algebra
and then $\Cal H^1$ is a nilpotent group.\fin

\begin{coro} Let $\Cal G$ a  solvable subgroup of $\dif{_1}$, then $\Cal
G^1=[\Cal G,\Cal G]$ is a nilpotent group.
\end{coro}

From the proposition \ref{expo} is clear that if $\Cal G\subset
\dih{_1}$ is solvable (nilpotent) if and only if the algebra
associate to $\Cal G$ by the $\exp $ function has to be solvable
(nilpotent). In particular in dimension 2 we have

\begin{prop}\label{nilp} Every nilpotent subalgebra $\Cal L$ of   $\camfdos{}$ is
metabelian.
\end{prop}
\proof Let $\Cal R$ be the center of $\Cal L$.  Since $\Cal R$ is
non trivial then $\Cal R\otimes \hat K(\C^2) $ is a vector space
of dimension 1 or 2 over $\hat K(\C^2) $, where $\hat K(\C^2) $ is
the fraction field of $\Cal O(\C^2)$.

In the case when the dimension is 2, there exist formal fields
$\mathfrak f$ and $\mathfrak g$ linearly independent over $\hat
K(\C^2) $, thus every element $\mathfrak h\in \Cal L$, can be
written as $\mathfrak h=u\mathfrak f+v\mathfrak g$, and since
$\mathfrak f,\mathfrak g\in \Cal R$ follows $\mathfrak
f(u)=\mathfrak f(v)=\mathfrak g(u)=\mathfrak g(v)=0$, i.e. $u$ and
$v$ are constants, therefore $\Cal L$ is abelian algebra.

If the dimension of the center is 1, let $\mathfrak f$ be a
non-trivial element of $\Cal R$, and $\Cal S=\hat K(\C^2)\mathfrak
f\cap \Cal L$, is clear that $\Cal R\subset \Cal S$ is an abelian
subalgebra of $\Cal L$. In the case $\Cal L=\Cal S$ we have
nothing to proof. Otherwise, since $\Cal S$ is an ideal  of $\Cal
L$, let $\mathfrak g$ be an element of $\Cal L$ such that its
image at $\Cal L/\Cal S$ is in the center of $\Cal L/\Cal S$. In
the same way every element of $\Cal L$ is of the form $u\mathfrak
f+v\mathfrak g$ where $\mathfrak f(u)=\mathfrak f(v)=0$, moreover
since $[u\mathfrak  f+v\mathfrak g,\mathfrak f]\in\Cal S$, i.e.
$\mathfrak g(v)=0$, follows that $v$ is constant, and then $\Cal
L$ is a metabelian algebra.\fin

\begin{coro}
Let $\Cal G$ a subgroup of $\dihdos{}$, then
\begin{enumerate}
\item[a)] If $\Cal G$ is solvable then $\Cal G^7=\{id\}$.
\item[b)] If $\Cal G$ is nilpotent then $\Cal G^6=\{id\}$.
\end{enumerate}
\end{coro}

Observe that $b)$ is  weaker that the Ghys theorem

\begin{teo}[Ghys] Let $\Cal G$ be a nilpotent subgroup of $\difdos{}$, then
$\Cal G$ is a metabelian group.
\end{teo}

The proposition \ref{nilp} gives a characterization of abelian
subgroup of $\dihdos{_1}$.

\begin{coro} If $\Cal G<\dihdos{_1}$ be a abelian group, then one
of the following items are true
\begin{enumerate}
\item $\Cal G< \left\langle \exp(N\Cal X)\in \dihdos{_1} \left|{\hbox {$N$ rational holomorphic
function} \atop \hbox {such that  $\Cal X
(N)=0$}}\right.\right\rangle$ where $ \exp(\Cal X)\in \Cal G$.
\item
$\Cal G< \langle F^{[t]}\comp G^{[s]}|t,s\in \C\rangle$, where
$F,G\in \Cal G$ and $[F,G]=Id$.
\end{enumerate}
\end{coro}

If  $F$ and $G$ are two elements of an abelian group of
diffeomorphisms tangent to the identity then from the lemma
\ref{comuta} we have
\begin{equation}
\nabla F_{r+1}(X){\cdot}G_{s+1}(X) -\nabla G_{s+1}(X)\cdot
F_{r+1}(X)\equiv 0. \label{comuta3}
\end{equation}
It is clear that in the particular case when the dimension is  1
the equation (\ref{comuta3}) is equivalent to say $r=s$. In
general, this is false in dimension $>1$, for example
$F(X)=\exp(\mathfrak f)(X)\in \dihdos{_2}$ and
$G(X)=\exp(\mathfrak g)(X)\in\dihdos{_3}$, where $$\mathfrak f=
(x^2+3xy)\parc x +(3xy+y^2)\parc y$$ and $$\mathfrak
g=(3x^3-5x^2y+xy^2+y^3)\parc x+( x^3+x^2y-5xy^2+3y^3)\parc y.$$
Since $[\mathfrak f,\mathfrak g]=0$ it follows that $F(X)$ and
$G(X)$ commute. In addition, $F$ and $G$ are holomorphic
diffeomorphisms because $\mathfrak f$ and $\mathfrak g$ are
holomorphic.

\begin{defi}
Let $F\in \dih{_{r+1}}$. $F$ is called {\bf dicritic}  if
$F(X)=X+F_{r+1}(X)+\cdots$, where $F_{r+1}(X)=f(X) X$ and $f$ is a
homogeneous polynomial of degree $r$.
\end{defi}

We are going to prove   a generalization of the dimension one
classification of abelian group for dicritic diffeomorphisms.

\begin{prop}\label{pr:comuta}
  Let $F\in \dif{_{r+1}}$ and $G\in \dif{_{s+1}}$.  Suppose that  $F$ is a
dicritic diffeomorphism, and $F(G(X))=G(F(X))$, then $r=s$ and $G$
is also a dicritic diffeomorphism.
\end{prop}
\proof For the lemma \ref{comuta} we have
\begin{equation}\eqalign
{\nabla F_{r+1}(X) G_{s+1}(X) -&\nabla G_{s+1}(X) F_{r+1}(X)=\cr
&=(f(X)I+(x_i\pa f{x_j})){\cdot}G_{s+1}(X) -\nabla G_{s+1}(X)\cdot
f(X)X\cr &=(f(X)I+(x_i\pa f{x_j})){\cdot}G_{s+1}(X)
-(s+1)f(X)G_{s+1}(X)\cr &=(-sf(X)I+(x_i\pa f{x_j}))\cdot
G_{s+1}(X)\label{comuta4}}\end{equation} From the identity
$\det(aI+AB)=\det(aI+BA)$

$$\eqalign{ \det(&-sf(X)I+(x_i\pa f{x_j}))\cr &=
\left|\matrix{-sf(X)+x_1\pa f{x_1}&x_2\pa f{x_2}&\cdots&x_n\pa
f{x_n}\cr x_1\pa f{x_1}&-sf(X)+x_2\pa f{x_2}&\cdots&x_n\pa
f{x_n}\cr \vdots&\vdots&\ddots&\vdots\cr x_1\pa f{x_1}&x_2\pa
f{x_2}&\cdots&-sf(X)+x_n\pa f{x_n}\cr}\right|\cr &=
\left|\matrix{-sf(X)&0&\cdots&sf(X)\cr 0&-sf(X)&\cdots&sf(X)\cr
\vdots&\vdots&\ddots&\vdots\cr x_1\pa f{x_1}&x_2\pa
f{x_2}&\cdots&-sf(X)+x_n\pa f{x_n}\cr}\right|\cr &=
\left|\matrix{-sf(X)&0&\cdots&0\cr 0&-sf(X)&\cdots&0\cr
\vdots&\vdots&\ddots&\vdots\cr x_1\pa f{x_1}&x_2\pa
f{x_2}&\cdots&(r-s)f(X)\cr}\right|\cr &=(-s)^{n-1}(r-s)(f(X))^n}
$$ So multiplying by the adjoint  matrix of $-sf(X)I+(x_i\pa
f{x_j})$ we obtain $(-s)^{n-1}(r-s)(f(X))^nG_{s+1}=0$ and  then
$r=s$. Denoting $G_{r+1}=(g_1,\dots, g_n)^t$, it is easy to see,
multiplying the $i$ row in (\ref{comuta4}) by $x_n$ and
subtracting the last row multiplied by $x_i$, that
$-sf(X)x_ng_i+sf(X)x_ig_n=0$. Therefore $x_n|g_n$.  Defining
$g=\frac {g_n}{x_n}$, we conclude $g_i=x_ig$ for all $i$.\fin

\begin{teo}\label{5.5}
 Let $\Cal G<\dihdos{_1}$ be abelian group, and $F\in \Cal G$ dicritic
diffeomorphism. Suppose that $\exp(\mathfrak f)(x,y)=F(x,y)$ where
$$\mathfrak f=(f(x,y)x+p_{k+2}(x,y)+\cdots)\parc x +
(f(x,y)y+q_{k+2}(x,y)+\cdots)\parc y,$$ $f(x,y)$ is a homogeneous
polynomial of degree $k$ and
$g.c.d(f,xq_{k+2}(x,y)-yp_{k+2}(x,y))=1$, then $$\Cal G<\langle
\exp(t\mathfrak f)(x,y)|t\in \C\rangle.$$
\end{teo}
\proof Let $G\in \Cal G$, for the proposition \ref{pr:comuta}, $G$
is a $k+1$-flat dicritic diffeomorphism. Let $\mathfrak g$ such
that $\exp(\mathfrak g)(x,y)=G$, then $$\mathfrak g=g(X)\vec
R+(s_{k+2}+\cdots)\parc x+ (t_{k+2}+\cdots)\parc y, $$ where $\vec
R=x\parc x+y\parc y$ and $X=(x,y)$. Since $[\mathfrak f,\mathfrak
g]=0$, each $j$-jet has to be zero. The $2k+2$-jet is $[f\vec
R,s_{k+2}\parc x+ t_{k+2}\parc y]+[p_{k+2}\parc x+ q_{k+2}\parc
y,g\vec R]=0$,  by a straightforward calculation we have
$$[f\vec R,s_{k+2}\parc x+ t_{k+2}\parc
y]=(k+1)f(s_{k+2}\parc x+ t_{k+2}\parc y)-(s_{k+2}\pa fx+
t_{k+2}\pa fy)\vec R,$$ therefore
$$\eqalign{(k+1)(fs_{k+2}-gp_{k+2})&=x(s_{k+2}\pa fx+ t_{k+2}\pa
fy-p_{k+2}\pa gx- q_{k+2}\pa gy)\cr
(k+1)(ft_{k+2}-gq_{k+2})&=y(s_{k+2}\pa fx+ t_{k+2}\pa
fy-p_{k+2}\pa gx- q_{k+2}\pa gy)}$$ In particular, $$\frac
{fs_{k+2}-gp_{k+2}}x=\frac{ft_{k+2}-gq_{k+2}}y,$$ or equivalently,
$f(xs_{k+2}-yt_{k+2})=g(xq_{k+2}-yp_{k+2})$, but
$gcd(f,xq_{k+2}-yp_{k+2})=1$, it follows that $f|g$, and since $f$
and $g$ have  the same degree then $g=rf$ where $r\in \C$.
Substituting $g$ we obtain the system of equations $$
\pmatrix{(k+1)f-x\pa fx&-x\pa  fy\cr -y\pa fx&(k+1)f-y\pa
fy}\pmatrix{s_{k+2}-rp_{k+2}\cr t_{k+2}-rq_{k+2}}=\pmatrix{0\cr
0}$$ multiply by the adjoint matrix
$(k+1)f^2\pmatrix{s_{k+2}-rp_{k+2}\cr
t_{k+2}-rq_{k+2}}=\pmatrix{0\cr 0}$ then $s_{k+2}=rp_{k+2}$ and
$t_{k+2}=rq_{k+2}$. Notice  that this last calculation is true in
arbitrary dimension. Finally,  suppose that $s_{k+j}=rp_{k+j}$ and
$t_{k+j}=rq_{k+j}$ for $j=2,\dots,i$. The $(2k+i+1)$-jet of the
bracket is $[f\vec R,s_{k+i+1}\parc x+ t_{k+i+1}\parc
y]+[p_{k+i+1}\parc x+ q_{k+i+1}\parc y,g\vec R]=0$ because the
symmetrical terms of the summation $$\sum_{j=2}^i[p_{k+j}\parc x+
q_{k+j}\parc y,s_{k+i+2-j}\parc x+ t_{k+i+2-j}\parc y]$$ is zero.
Then the equalities $s_{k+i+1}=rp_{k+i+1}$ and
$t_{k+i+1}=rq_{k+i+1}$ follows in the  same way to the case
$k+2$.\fin

Observation: The condition over $F$ is generic and means that
$(0,0)\in \C^2$ is an isolated singularity of $ \mathfrak
f_0=(f(x,y)x+p_{k+2}(x,y))\parc x + (f(x,y)y+q_{k+2}(x,y))\parc
y$. With a similar  condition the theorem is true in arbitrary
dimension.

\section{Convergent Orbits}

\begin{defi} Let $F\in \dih{_{k+1}}$, where $F(X)=X+F_{k+1}(X)+\cdots$.
 We say that $[V]\in \C P(n-1)$ is a characteristic direction of
$F$,  if there exists $\lambda\in \C$ such that
$$F_{k+1}(V)=\lambda V$$ Moreover  the direction $[V]$ is called
non degenerate if  $\lambda\ne 0$.
\end{defi}

\begin{defi} Let $F\in \dih{_2}$. A parabolic curve for $F$ at the origin is
an injective  map $\varphi:D_1\to \C^n$, where $D_1=\{x| \
|x-1|\le 1\}$, holomorphic in $int(D_1)$, such that
$\varphi(0)=0$, $\varphi(D_1)$ is invariant under $F$ and
$(F|_{D_1})^n\to 0$ when $n\to \infty$
\end{defi}

\begin{teo}[Hakim \cite{Ha}]
Let $F\in \dih{_{k+1}}$, then for every non-degenerate
characteristic direction $[V]$ there are a  parabolic curves
tangent to $[V]$ at the origin.
\end{teo}

\begin{teo}[Abate \cite{Ab}]
Let $F\in \dihdos{_{k+1}}$ such that the origin is a isolated
fixed point. Then there exist $k$ parabolic curves for $F$ at the
origin.
\end{teo}

Denote $r(v)=q_{k+1}(1,v)-vp_{k+1}(1,v)$ and $p(v)=p_{k+1}(1,v)$.

\begin{lem}\label{seq1} Let $F\in \dihdos{_{k+1}}$ and suppose that
 $(x_n,y_n)=F(x_{n-1},y_{n-1})$ is  a sequence converging to 0 such that
$\frac{y_n}{x_n}\to v$ when $n\to \infty$.  Then $$\lim_{n\to
\infty} \frac1{nx_n^k}=-kp(v)$$
\end{lem}

\proof After a blow up at $0\in \C^2$ in the chart $y=xv$, we
obtain the diffeomorphism $$\vect{x_1}{v_1}=\tilde
F(x,v)=\vect{x+x^{k+1}p(v)+x^{k+2}p_{k+2}(1,v)+\cdots}
{v+x^kr(v)+\cdots}$$ Rewriting the first equation $$\eqalign{
\frac 1{x_1^k}&=\frac
1{(x+x^{k+1}p(v)+x^{k+2}p_{k+2}(1,v)+\cdots)^k}\cr &=\frac
1{x^k}(1+x^{k}p(v)+x^{k+1}p_{k+2}(1,v)+\cdots)^{-k}\cr &=\frac
1{x^k}-kp(v)+o(x). } $$ Let's define $(x_j,
v_j)=F(x_{j-1},v_{j-1})$. From the equation above we get the
telescopic sum $$ \frac 1{x_n^k}-\frac 1{x^k}=\sum_{j=1}^n \frac
1{x_j^k}-\frac 1{x_{j-1}^k} =-\sum_{j=1}^n
(kp(v_{j-1})+o(x_{j-1}))$$ Divide by $n$ and make $n$ tends to
$\infty$, we deduce $$\lim_{n\to \infty}
\frac1{nx_n^k}=-k\lim_{n\to \infty}\frac 1n\sum_{j=1}^n
kP(v_{j-1}) +\lim_{n\to \infty}\sum_{j=1}^n \frac 1n
o(x_{j-1})=-kp(v).\eqno \fin$$

\begin{prop}\label{R(v)=0} Let $(x_0,y_0)\in U$ open neighborhood of 0 such
that the sequence $(x_j, y_j)=F(x_{j-1},y_{j-1})$ converge to 0,
and $\frac {y_j}{x_j}$  converge to $v\in\C P(1)=\overline \C$,
then $r(v)=0$. \end{prop} \proof From the lemma \ref{seq1}. we
have $\lim_{n\to \infty} \frac1{nx_n^k}=-kp_{k+1}(1,v)$ and the
same way $\lim_{n\to \infty} \frac1{ny_n^k}=-kq_{k+1}(\frac
1v,1)$. Dividing these relations we get $$\frac 1{v^k}=\lim_{n\to
\infty} \frac{x_n^k}{y_n^k}=\frac {\lim_{n\to \infty}
\frac1{ny_n^k}}{\lim_{n\to \infty}
\frac1{nx_n^k}}=\frac{q_{k+1}(\frac 1v,1)}{p_{k+1}(1,v)} ,$$
therefore $v^{k+1}q_{k+1}(\frac 1v,1)-vp_{k+1}(1,v)=r(v)=0$.\fin

\begin{teo}\label{7.3} Let  $F:(\C^2,0)\to (\C^2,0)$ be a dicritic diffeomorphism fixing zero,
 i.e. $F$ can be represented by a convergent series
$$F(x,y)=\left(\matrix{x+xp_k(x,y)+p_{k+2}(x,y)+\cdots\cr
 y+yp_k(x,y)+q_{k+2}(x,y)+\cdots}\right),$$
and $\tilde F=\Pi^*F:(\tilde\C^2,D)\to (\tilde\C^2,D)$ the
continuous extension of the diffeomorphism  after  making the
blow-up in $(0,0)$. Then there exist open sets $U^+,U^-\subset
\tilde \C^2$ such that
\begin{enumerate}
\item $\overline{U^+\cup U^-}$ is a neighborhood of
$D\setminus\{(1:v)\in D|p_k(1,v)=0\}$. \item For all $P\in U^+$,
the sequence  $\{\tilde F^n(P)\}_{n\in \N}$ converge and
$\lim\limits_{n\to\infty}\tilde F^n(P)\in D$. \item For all $P\in
U^-$, the sequence  $\{\tilde F^{-n}(P)\}_{n\in \N}$ converge and
$\lim\limits_{n\to\infty}\tilde F^{-n}(P)\in D$.
\end{enumerate}
\end{teo}

\proof Making a blow  up  at $(0,0)$, and regarding the
diffeomorphism in the chart $(x,v)$,  $v=\frac yx$ we obtain
$${x_1\choose v_1}=\tilde
F(x,v)=\left(\matrix{x+x^{k+1}p(v)+\sum_{j=2}^\infty
a_{k+j}(v)x^{k+j} \cr
 v+\sum_{j=1}^\infty b_{k+j}(v)x^{k+j}}\right),$$
where $a_j$ and $b_j$ are polynomial of degree  less than $j+2$.
Let $q$ be a arbitrary point in $D\setminus\{(1:v)\in D|p(v)=0\}$.
We can suppose making a linear change of coordinates that
$q=(0,0)\in \tilde \C^2$. Since $\tilde F$ is holomorphic in some
neighborhood of $(0,0)$ there exist $r_1,r_2>0$ and constants
$C_1,C_2>0$ such that $$ ||\sup_{|v|\le r_1}a_j(v)||\le C_1
r_2^j,\hbox{ and } ||\sup_{|v|\le r_1}b_j(v)||\le C_2 r_2^j\hbox{
for all $j$}.$$ Now making a ramificated change of coordinates
$w=\frac 1{x^k}$, it follows that $${w_1\choose v_1}=\overline
F(w,v)= \left(\matrix{w-kp(v)+\sum_{j=1}^\infty c_j(v)\frac 1
{w^{\frac jk}} \cr
 v+\sum_{j=1}^\infty b_{k+j}(v)\frac 1{w^{1+\frac jk}}}\right),$$
is holomorphic in $\{(w,v)\in \C_k{\times}\C|\ |v|< r_1,\
|w|>r_2^k\}$, where $\C_k$ is a $k$-fold covering of $\C^*$. In
particular there exist $C_3,C_4  >0$ such that $$||\sup_{|v|\le
r_1}c_j(v)||\le C_3 r_2^j \hbox{ for all $j$ and}
\left|\sum_{j=1}^\infty b_{k+j}(v)\frac 1{w^{1+\frac
jk}}\right|<\frac {C_4}{|w|^{1+\frac 1k}} .$$ Choose $r_1>r>0$
such that $\left|\frac {p(v)}{p(0)}-1\right|<\frac 14$ for all
$|v|<r$, and $$ R>\max \left\{r_2^k\left(\frac
{4C_3}{k|p(0)|}+1\right)^k, \frac {(2C_4)^k}{r^k} \left(1+\frac
{4C(k)}{k|p(0)|}\right)^k,1\right\}, $$ where
$C(k)=\displaystyle\int_0^{\infty} \frac 1{(1+x^2)^{\frac
{k+1}{2k}}}d x<k+1$. Let $$V_q^+(a)= \left\{(w,v)\in \C_k\times
\C\Bigl|\ |v|<a,\quad \left|\arg\left(-\frac
w{p(0)}-\frac{2R}{|p(0)|}\right)\right|<\frac {2\pi}3\right\}, $$
where $\arg$ is defined  from $\C_k^*$ to $(-\pi,\pi]$. It is easy
to prove that $V_q^+(r)\subset \{(w,v)\in \C_k{\times}\C|\ |v|< r,\
|w|>R\}$.

\

Claim: For all $p=(w_0,v_0)\in V_q^+(\frac r2)$ we have
\begin{enumerate}
\item $\overline F^n(p)=(w_n,v_n)\in V_q^+(r)$ for all $n\in \N$.
\item The sequence $\{v_n\}_{n\in \N}$ converge. \item $|w_n|\to
\infty$ when $n\to +\infty$.
\end{enumerate}

{\it Proof of the claim:} For a arbitrary point  $(w,v)$ in
$V_q^+(r)$ we have $$ \eqalign{|w_1-(w-kp(0))|&\le
\left|k(p(0)-p(v))+\sum_{j=1}^\infty c_j(v)\frac 1 {w^{\frac
jk}}\right|\cr &\le k|p(0)-p(v)|+\sum_{j=1}^\infty C_3\frac
{r_2^j} {|w|^{\frac jk}}<k\frac {|p(0)|}4 +C_3\frac {r_2}
{|w|^{\frac 1k}-r_2}\cr &<k\frac {|p(0)|}4+C_3\frac
{r_2}{r_2\left(\frac {4C_3}{k|p(0)|}+1\right)-r_2}=k\frac
{|p(0)|}2.} $$ Let $p=(w_0,v_0)\in V_q^+(\frac r2)$ and suppose
that $(w_i,v_i)\in V_q^+(r)$ for \hbox{$i=0,\dots,j-1$}. Then
$$|w_j-(w_0-jkp(0))|=\Bigl|\sum_{i=1}^j
w_i-(w_{i-1}-kp(0))\Bigr|\le \sum_{i=1}^j |w_i-(w_{i-1}-kp(0))|<jk
\frac {|p(0)|}2$$ therefore  $w_j=w_0-jkp(0)(1+\delta_j)$ where
$|\delta_j|<\frac 12$, and $$-\frac
{w_j}{p(0)}-\frac{2R}{|p(0)|}=-\frac
{w_0}{p(0)}-\frac{2R}{|p(0)|}+jk(1+\delta_j),$$ but since
$|\arg(1+\delta_j)|<\frac \pi3$ it follows that $|\arg(-\frac
{w_j}{p(0)}-\frac{2R}{|p(0)|})|<\frac {2\pi}3$. Now notice that
$|v_i|<|v_{i-1}|+\frac {C_4}{|w_{i-1}|^{1+\frac 1k}}$ for all
$i=1,\dots, j$, therefore $|v_j|<|v_0|+\sum_{i=1}^j \frac
{C_4}{|w_{i-1}|^{1+\frac 1k}},$ but
\begin{equation}
\eqalign{\sum_{i=1}^j \frac 1{|w_{i-1}|^{1+\frac 1k}}& <
\sum_{i=-\infty}^\infty \frac 1{(R^2+(\frac {ik|p(0)|}2)^2)^{\frac
{k+1}{2k}}}\cr &< \frac 1{R^{1+\frac 1k}}+2\int_0^\infty \frac
1{(R^2+(\frac {k|p(0)|}2)^2x^2)^{\frac {k+1}{2k}}}dx\cr &<\frac
1{R^{\frac 1k}}+\frac 4{k|p(0)| R^{\frac 1k}}\int_0^{\infty} \frac
1{(1+x^2)^{\frac {k+1}{2k}}}d x=\frac 1{R^{\frac 1k}}\left(1+\frac
{4C(k)}{k|p(0)|}\right)\cr &<\frac r{2C_4}, }\label{desig1}
\end{equation}
it follows that $|v_j|<r$, and therefore $(w_j,v_j)\in V_q^+(r)$.
Now we have trivially that $|w_j|=| w_0-jkp(0)(1+\delta_j)|\to
\infty$ when $n\to\infty$ and $\{v_j\}_{j\in \N}$ is a Cauchy
sequence because $\sum_{i=0}^\infty\frac {C_4}{|w_i|^{1+\frac
1k}}$ is a Cauchy series.

In the same way we can obtain the open set $V_q^-$  changing  in
the proof the diffeomorphism germ $F$ by $F^{-1}$. Finally we
conclude the proof making $U^+=\bigcup V_q^+$ and $U^-=\bigcup
V_q^-$ where $q\in D\setminus\{(1:v)\in D|p_k(1,v)=0\}$.\fin

\begin{teo} Let $F\in \dihdos{_{k+1}}$, where
$F(x,y)=\left(\matrix{x+p_{k+1}(x,y)+\cdots\cr
 y+q_{k+1}(x,y)+\cdots}\right)$ and
 $r(v)=vp_{k+1}(1,v)-q_{k+1}(1,v)\not\equiv 0$.
 Suppose that $v_0\in \C$ satisfies $r(v_0)=0$ and
$\Re(\frac {r'(v_0)}{p(v_0)})>0$. Then there exist open sets $U^+$
and $U^-$ such that $(0,0)\in \partial U^{\pm}$ and for each point
$(a^+,b^+)\in U^+$ and $(a^-,b^-)\in U^-$ the sequences
$(a^\pm_n,b^\pm_n)=F^{{\pm}n}((a^\pm,b^\pm))$ converge and
$$\lim_{n\to \infty} F^{{\pm}n}(a^\pm,b^\pm)= 0\hbox{ and }
\lim_{n\to \infty} \frac {b_n^+}{a_n^+}=\lim_{n\to \infty} \frac
{b_n^-}{a_n^-}=v_0.$$
\end{teo}

\proof For a linear change of coordinates we can suppose that
$v_0=0$.
Let $r(v)=vs(v)$, $\alpha=s(0)$,  $\beta=p(0)$ and $\tilde
F:(\tilde\C^2,D)\to (\tilde\C^2,D)$ the continuous extension of
the diffeomorphism  after  making the blow-up $\Pi$ at  $(0,0)$,
i.e. in the chart $(x,v)$ where $y=vx$ $$(x_1,v_1)=\tilde
F(x,v)=\left(\matrix{x+x^{k+1}p(v)+x^{k+2}(\cdots)\cr
 v+x^kr(v)+x^{k+1}(\cdots)}\right).$$
Since $\tilde F|_D=id|_D $  we can forget the dynamic in
$D=\{x=0\}$ and make a ramificate change of coordinates $w=\frac
1{x^k}$. Then $F$ is representing in the $(w,v)$ coordinates
system as $$(w_1,v_1)=\overline F(x,v)=\left(\matrix{w-kp(v)+
o(\frac 1{w^{ 1/k}})\cr
 v+\frac {s(v)}w v+ o(\frac 1{w^{1+ 1/k}})}\right).$$
In particular, there exists a neighborhood $A=\{(w,v)\in
\C_k^*{\times}\C| \ |w|>R_1, \ |v|<r_1\}$ of $(\infty,0)$
 and constants $C_1,C_2>0$
such that $$|w_1-(w-kp(v))|<C_1 \frac 1{|w|^{1/k}}\ \hbox{ and }\
\left|v_1-v\Bigl(1+ \frac{s(v)}w\Bigr)\right|<C_2 \frac
1{|w|^{1+1/k}}.$$ Let denote $\theta=\frac \pi2-|\arg(\frac
\alpha\beta)|>0$. Choose $0<r<r_1$ such that $$\left|\frac
{p(v)}\beta-1\right|<\frac 12 \sin \frac \theta 2,\quad \left|\arg
\frac{s(v)}\alpha\right|<\frac \theta 2$$ and  $R>\max
\left\{R_1+1,\left( \frac {2C_1}{|\beta| k \sin \frac \theta 2}
\right)^k, \frac{(2C_2)^k}{r^k} \left(1+\frac {2C(k)}{k(1- \sin
\frac \theta 2)|\beta|}\right)^k\right\}.$ Let define
$$V^+(a)=\{(w,v)\in A| \ |v|<a, \  \Re (\frac w\alpha e^{\pm\frac
\theta2 i}) <-R\}.$$

See that for all $(w,v)\in V^+(r)$ we have $\Re (\frac
w{s(v)})<-R$, i.e. $\frac {s(v)}w$ is in the disc of center
$-\frac 1R$ and radius  $\frac 1R$, it follows that
$$\left|1+\frac {s(v)}w\right|<1.$$ In addition $$
\eqalign{|w_1-(w-k\beta)|&\le |k(\beta-p(v))|+C_1 \frac
1{|w|^{1/k}}\cr &\le \frac {k|\beta|}2 \sin \frac \theta 2+\frac
{k|\beta|}2 \sin\frac\theta2\cr &=k|\beta| \sin \frac \theta 2.}
$$ Let $p=(w_0,v_0)\in V^+(\frac r2)$ and suppose that
$(w_i,v_i)\in V^+(r)$ for \hbox{$i=0,\dots,j-1$}. Then
$$|w_j-(w_0-jk\beta)|\le \sum_{i=1}^j |w_i-(w_{i-1}-k\beta)|
<jk|\beta |\sin  \frac \theta2,$$ therefore
$w_j=w_0-jk\beta(1+\delta_j)$ where $|\delta_j|<\sin \frac
\theta2$. Now since $\frac {w_j}{\alpha}=\frac
{w_0}{\alpha}-jk\frac \beta\alpha(1+\delta_j),$ and $$\left|\arg
(\frac \beta\alpha(1+\delta_j))\right|\le \left|\arg \frac
\beta\alpha\right|+|\arg (1+\delta_j)|\le \frac \pi2-\theta +\frac
\theta2 =\frac \pi2-\frac \theta2,$$ it follows that $\Re(\frac
\beta\alpha(1+\delta_j)e^{{\pm}\frac \theta 2 i}) \ge 0 $ and
therefore $\Re(\frac {w_j}{\alpha} e^{{\pm}\frac \theta 2 i})<-R$.

In order to bound $v_j$, see that $$|v_i|<\left|\left(1+\frac
{s(v_{i-1})}{w_{i-1}}\right)v_{i-1}\right| +\frac
{C_2}{|w_{i-1}|^{1+\frac 1k}}\le |v_{i-1}| +\frac
{C_2}{|w_{i-1}|^{1+\frac 1k}}$$ for all $i=1,\dots, j$, using a
similar calculation as  the inequality in (\ref{desig1}) we have
$$|v_j|<|v_0|+\sum_{i=1}^j \frac {C_2}{|w_{i-1}|^{1+\frac
1k}}<r.$$

It follows that $(w_j,v_j)\in V^+(r)$, $|w_j|=|
w_0-jkp(0)(1+\delta_j)|\to \infty$ when  $n\to\infty$ and
$\{v_j\}_{j\in \N}$ is a Cauchy sequence because
$\sum_{i=0}^\infty\frac {C_2}{|w_i|^{1+\frac 1k}}$ is a Cauchy
series. Moreover from the proposition \ref{R(v)=0} we know that
$\{v_j\}_{j\in \N}$ converge to 0. \fin

\section{Normal forms}

Let $\tilde F,\tilde G\in \dihdos{_{k+1}}$ such that their
$(k+1)$-jet are equal and suppose that they are formally
conjugate. Let $F(x,v)=\left(\matrix{x+x^{k+1}p(v)+\cdots\cr
v+x^kr(v)+\cdots}\right)$ and
$G(x,v)=\left(\matrix{x+x^{k+1}p(v)+\cdots\cr
v+x^kr(v)+\cdots}\right)$ be the diffeomorphisms obtained after
the blow-up at 0 where $r(v)$ and $p(v)$ are polynomials of degree
$(k+2)$ and  $(k+1)$ respectively. Let $H$ be
 a formal  diffeomorphism  that conjugates $F$ and $G$. Notice that
$H(x,v)=\left(\matrix{x+\sum_{l=1}^\infty x^{l+1}h_{1,l}(v)\cr
v+\sum_{l=1}^\infty x^l h_{2,l}(v)}\right)$ where $h_{1,l}$ and
$h_{2,l}$ are polynomial of degree $l+1$ and $l+2$ respectively,
in particular $h_{1,l}$ and $h_{2,l}$ are holomorphic functions in
$\C$. The following definition is borrowed from \cite{Vo1}.

\begin{defi} A formal Taylor series in $x$ is called semiformal in $U$ if its
coefficients holomorphically depend on $v$ in the same domain $U$.
A formal change of coordinates is called  semiformal in $U$ if its
components are semiformal series in $U$.
\end{defi}

It is clear that if $H$ is a semiformal change of coordinates in
$\overline\C$, then $F$ and $G$ are semiformally conjugates in
$\overline \C$. We are interested in finding some semiformal
invariant, for that we need the following lemma.

\begin{lem}\label{eq:difer} Let $F$ and $G$ as above. Suppose that
$G-F=\left(\matrix{x^{k+l+1}\phi_1(v)+\cdots\cr
x^{k+l}\phi_2(v)+\cdots}\right)$ and that
$H(x,v)=\left(\matrix{x+x^{l+1}h_1(v) +\cdots\cr v+x^l
h_2(v)+\cdots}\right)$ conjugates $F$ and $G$, i.e. $G\comp
H=H\comp F$. Then
\begin{equation}
r(v)\pmatrix{h_1(v)\cr  h_2(v)}'- \pmatrix{(k-l)p(v)&p'(v) \cr
kr(v)&r'(v) -lp(v) } \pmatrix{h_1(v)\cr
h_2(v)}=\pmatrix{\phi_1(v)\cr  \phi_2(v)}.\label{eq:difer1}
\end{equation}
\end{lem}

\proof
 For the lemma \ref{comuta} we have
{\small $$\eqalign{ 0&=G\comp H-H\comp F =G\comp H-H\comp
(G+(F-G))\cr &=G\comp H-H\comp G -\nabla H(G)\cdot
(F-G)+\left({x^{2k+2l+2}(\cdots)\atop x^{2k+2l}(\cdots)}\right)\cr
&=\pmatrix{(k+1)x^kp(v)&x^{k+1}p'(v)\cr kx^{k-1}r(v)&x^kr'(v)}
\pmatrix{x^{l+1}h_1(v)\cr x^l h_2(v)}
-\pmatrix{(l+1)x^lh_1(v)&x^{l+1}h_1'(v)\cr
lx^{l-1}h_2(v)&x^lh_2'(v)} \pmatrix{x^{k+1}p(v)\cr x^k r(v)}\cr
&\hskip 8cm +\left({x^{k+l+1}\phi_1(v)\atop
x^{k+l}\phi_2(v)}\right) +\left({x^{k+l+2}(\cdots)\atop
x^{k+l+1}(\cdots)}\right)\cr
&=\pmatrix{x^{k+l+1}((k-l)p(v)h_1(v)+p'(v)h_2(v)
-r(v)h_1'(v)+\phi_1(v))\cr x^{k+l}(kr(v)h_1(v)+(r'(v)
-lp(v))h_2(v)-r(v)h_2'(v)+\phi_2(v)) }
+\left({x^{k+l+2}(\cdots)\atop x^{k+l+1}(\cdots)}\right) } $$} In
particular, we obtain the linear differential equation
(\ref{eq:difer1}).\fin

\begin{prop}Let $F$ and $G$  as in
the lemma \ref{eq:difer}. Suppose that $r(v)\not\equiv 0$ and  $U$
is a simply connected open set of $\C\setminus\{r(v)=0\}$. Then
(\ref{eq:difer1}) has holomorphic solution in $U$.
\end{prop}

\proof It follows from the fundamental existence theorem for the
holomorphic differential equations (See \cite{Wa} Theorem 2.1).
\fin

\begin{prop} Let $v_0$ be a root of $r(v)$. Suppose that $ \frac
{p(v_0)}{r'(v_0)}\notin \Q$. Then there exists a  holomorphic
solution $(h_1,h_2)$  of (\ref{eq:difer1}) in some neighborhood of
$v_0$ if and only if one of the following condition is true
\begin{enumerate}
\item $l\ne k$.
\item $l=k$ and $p'(v_0)\phi_2(v_0)=(r'(v_0)-kp(v_0))\phi_1(v_0)$.
\end{enumerate}
Moreover, that neighborhood is independent of the integer $l$.
\end{prop}

\proof Let $r(v)=(v-v_0)s(v)$ and $a=\frac {p(v_0)}{s(v_0)}\notin
\Q$, we are interested in   finding a solution of
\begin{equation}
(v-v_0)X'-\pmatrix{(k-l)\frac{p(v)}{s(v)}&\frac {p'(v)}{s(v)} \cr
k(v-v_0)&\frac{r'(v) -lp(v)}{s(v)} } X =\frac
1{s(v)}\pmatrix{\phi_1(v)\cr  \phi_2(v)}\label{solu1}
\end{equation}
in some neighborhood of $v_0$.

Since the difference between the eigenvalues of
$\pmatrix{(k-l)\frac{p(v_0)}{s(v_0)}&\frac {p'(v_0)}{s(v_0)} \cr
0&\frac{r'(v_0) -lp(v_0)}{s(v_0)} }$ is not a integer, it is known
(See \cite{Wa}) that there  exists a holomorphic change of
coordinates $Z=P(v)X$ ($P(v_0)=Id$) in some neighborhood
$V=B(v_0,R)$ of $v_0$    such that the system of differential
equations  (\ref{solu1}) is equivalent to
\begin{equation}
(v-v_0)Z'-\pmatrix{(k-l)a&a\frac {p'(v_0)}{p(v_0)} \cr 0&1-la } Z
=\pmatrix{\psi_1(v)\cr  \psi_2(v)}\label{solu2}
\end{equation}
In addition, $R$ only depend on the radius of convergence of the
series that defines locally
 each function in (\ref{solu1}) around $v_0$, and then independent of $l$.

Thus, we know that
 $\pmatrix{\psi_1(v)\cr
\psi_2(v)}=\sum\limits_{j=0}^\infty B_j (v-v_0)^j$ for all
$|v-v_0|<R$, and then  if we write $Z=\sum_{j=0}^\infty A_j
(v-v_0)^j$, we obtain for each $j$ a linear equation
$$\pmatrix{j-(k-l)a&-a\frac {p'(v_0)}{p(v_0)} \cr 0&j-(1
-la)}A_j=B_j.$$

Notice that, in the case $l\ne k$ or $j\ne 0$ this equation  has
solution because $a\notin \Q$, and when $l=k$ and $j=0$  we obtain
the linear equation $$\pmatrix{0&-\frac {p'(v_0)}{r'(v_0)} \cr
0&-1 +k\frac {p'(v_0)}{r'(v_0)}}A_j=\frac
1{r'(v_0)}\pmatrix{\phi_1(v_0)\cr \phi_2(v_0)},$$ which  has
solution if and only if $
p'(v_0)\phi_2(v_0)=(r'(v_0)-kp(v_0))\phi_1(v_0)$. Moreover
$$\lim_{j\to\infty}||A_j||^{\frac 1j}\le
\lim_{j\to\infty}\left\|\pmatrix{j-(k-l)a&-a\frac
{p'(v_0)}{p(v_0)} \cr 0&j-1 -la}^{-1}\right\|^{\frac 1j}
\cdot||B_j||^{\frac 1j}=\frac 1R,$$ Therefore  that the solution
is holomorphic in $B(v_0,R).$ \fin

\begin{teo}\label{semiconj}
Let $F(x,v)=\left(\matrix{x+x^{k+1}p(v)+\cdots\cr
v+x^kr(v)+\cdots}\right)$ and $v_0$ be a root of $r(v)$. Suppose
that $0\ne\frac {p(v_0)}{r'(v_0)}\notin \Q$. Then there exists
$\lambda_{v_0}\in \C$ such that $F$ is semiformally conjugate with
$$F_{\lambda_{v_0}}=\left(\matrix{x+x^{k+1}p(v)+\lambda_{v_0}x^{2k+1}\cr
v+x^kr(v)}\right)$$ in some neighborhood of $v_0$.
\end{teo}

\proof We  will go to construct the semiformal conjugation $H$ by
successive approximations. Define by induction the sequences
$F_l$ and $H_l$ from the initial condition  $F_0=F$ and $H_0=id$.

Now  for $l\ge 1$ and $l\ne k$, from the proposition
\ref{semiconj} we know that the system of differential equations
(\ref{eq:difer1}) in the lemma \ref{eq:difer}, where $G=F_{l-1}$
and
$$F_{l-1}-F_{\lambda_{v_0}}=\left(\matrix{x^{k+l+1}\phi_{1,l}(v)+\cdots\cr
x^{k+l}\phi_{2,l}(v)+\cdots}\right),$$ has holomorphic solution
$(h_{1,l},h_{2,l})$ in some neighborhood  $U$ that is independent
of $l$.

In the case $l=k$ see that
$$F_{k-1}-F_{\lambda_{v_0}}=\left(\matrix{x^{2k+1}(\phi_{1,k}(v)-
\lambda_{v_0})+\cdots\cr x^{2k}\phi_{2,k}(v)+\cdots}\right),$$
then  if we define $$\lambda_{v_0}=\phi_1(v_0)-\frac
{p'(v_0)\phi_2(v_0)}{r'(v_0)-kp(v_0)},$$ it follows for  the
proposition \ref{semiconj}, that there exists holomorphic solution
$(h_{1,k},h_{2,k})$ of (\ref{eq:difer1}).

In any case  we could define
$$H_l(x,v)=\left(\matrix{x+x^{l+1}h_{1,l}(v)\cr v+x^l
h_{2,l}(v)}\right)\hbox{ and } F_l=H_l^{-1}\circ F_{l-1}\circ
H_l.$$ See that in this case  $F_l$ and $H_l$ are holomorphic
diffeomorphisms in  $U$.

Moreover, $H=\lim_{n\to \infty}H_n\circ\cdots\circ H_0$ is a
semiformal diffeomorphism  that conjugates $F$ and
$F_{\lambda_{v_0}}$ in $U$. \fin

Since this conjugation is defined locally around some root of
$r(v)$, and we are interested in some global conjugation defined
at some neighborhood of the divisor, first we need to construct
some diffeomorphism that can be locally conjugate with $F$ around
every point of the divisor.

For that, let $L_F(v)$ denote the Lagrange interpolation
Polynomial of the points
$$(v_1,\lambda_1),\dots,(v_{k+2},\lambda_{k+2}),$$ where $v_j$ is
root of $r(v)=0$ and $\lambda_j$ is the constant found in the
theorem \ref{semiconj}. It is a simple consequence of the theorem
\ref{semiconj} that for all $U$ simply connected open set that
contains only one root of $r(v)=0$ there exists a semiformal
conjugation in $U$ that conjugates $F$ and
\begin{equation}
{x\choose v}\mapsto\left(\matrix{x+x^{k+1}p(v)+L_F(v)x^{2k+1}\cr
v+x^kr(v)}\right).\label{auxi}
\end{equation}

Notice that the diffeomorphism represented by (\ref{auxi}) does
not necessarily come from the blow up of some element of
$\dihdos{}$.

But it is easy to see, using the theorem \ref{semiconj} that the
blow up of $$F_L(x,y)=\left(\matrix{x+P_{k+1}(x,y)+L_F(\frac yx)
x^{2k+1}\cr y+Q_{k+1}(x,y)-x^{2k-1}yL_F(\frac yx)}\right)$$ is
locally semiformally  conjugate with (\ref{auxi}), and then it is
locally semiformally conjugate with $F$ for $k\ge 2$. Moreover, in
the chart $(x,v)$ the diffeomorphism  $F_L$ has the following
representation
$$\Pi^*F_L(x,v)=\left(\matrix{x+x^{k+1}p(v)+L_F(v)x^{2k+1}\cr
v+x^kr(v)-x^{2k}r(v)p(v)+x^{3k}(\cdots)}\right). $$

Let $\Cal U=\{U_i\}_{i=1,\dots,k+2}$ a   covering  of $\C
P(1)=\Pi^{-1}(0)$ such that $U_i$, $U_i\cap U_j$ are simply
connected  open sets,  such that no four of them have a nonempty
intersection and $v_i$ is in $U_j$ if and only if $i=j$. For each
$j$ there exists a semiformal diffeomorphism $H_j$ such that
$H_j\circ F=F_L\circ H_j$  in $U_j$.

Let $H_{ij}=H_j\circ H_i$ semiformal diffeomorphism defined in
$U_j\cap U_i$. Observe that each $H_{i,j}$ commutes with $F_L$.

\begin{teo} The cocycle $\{H_{ij}\}$ determine the class of formal conjugation
of $F$.
\end{teo}
First, we are going to prove that the cocycle $\{H_{ij}\}$
determine the class of semiformal conjugation of $F$. In fact
suppose that $\{H_{ij}=H'_j\circ {H'_i}^{-1}\}$ is the cocycle
associate to $G$. Then $$H_j\circ {H_i}^{-1}=H'_j\circ
{H'_i}^{-1}$$ in $U_i\cap U_j$, thus $H={H_j}^{-1}\circ
H_j'={H_i}^{-1}\circ H_i'$ is a global semiformal diffeomorphism
that conjugate $F$ and $G$. Let $$\vect xv\mapsto \vect
{x+\sum\limits_{j=k+1}^\infty a_{1j}(v)x^j}
{v+\sum\limits_{j=k}^\infty a_{2j}(v)x^j}$$ be the representation
of $H$ at the chart $(x,v)\in U_1$ and $$\vect ys\mapsto \vect
{y+\sum\limits_{j=k+1}^\infty b_{1j}(s)y^j}
{s+\sum\limits_{j=k}^\infty b_{2j}(s)y^j}$$ be the representation
of $H$ at the chart $(y,s)\in U_2$ where $U_1\simeq U_2\simeq
(\C,0){\times}\C$ and the change of coordinates between $U_1$ and
$U_2$ is given by $$\matrix{U_1\cap U_2&\to&U_1\cap U_2\cr
(x,v)&\mapsto&(y,s)=(vx,\frac 1v)}. $$ Therefore, we have that
$$\left(ys+\sum_{j=k+1}^\infty a_{1j}(\frac 1s)(ys)^j\right)
\left(\frac 1s+\sum_{j=k}^\infty a_{2j}(\frac
1s)(ys)^j\right)=y+\sum_{j=k+1}^\infty b_{1j}(s)y^j$$ and $$\frac
1{\frac1s+\sum_{j=k}^\infty a_{2j}(\frac 1s)(ys)^j}=
s+\sum_{j=k}^\infty b_{2j}(s)y^j$$ for every $(y,s)\in
(\C,0){\times}\C^*$.

Thus, it is easy to prove making the product and proceeding by
induction  that $a_{2j}(\frac 1s)s^{j-1}$ is holomorphic function
in $\C$ and therefore $a_{2j}(v)$ has to be a polynomial of degree
$\le j-1$, the same way $a_{1j}(v)$  has to be a polynomial of
degree $\le j$. Then, we can conclude that $H$ is a blow up of the
formal diffeomorphism in the variables $x,y$ that conjugate $F$
and $G$.\fin

Observe that in the  dicritic case,  i.e. $r\equiv 0$, the system
of  equations (\ref{eq:difer1}) of the lemma \ref{eq:difer}
reduces to
\begin{equation}
\pmatrix{(k-l)p(v)&p'(v) \cr 0 & -lp(v) } \pmatrix{h_1(v)\cr
h_2(v)}=\pmatrix{\phi_1(v)\cr  \phi_2(v)}\label{eq:dicrit}
\end{equation}

\begin{teo}\label{8.3} Let $\tilde F\in \difdos{_{k+1}}$ be dicritic diffeomorphism and
$F(x,v)=\left(\matrix{x+x^{k+1}p(v)+x^{k+2}(\cdots)\cr
v+x^{k+1}(\cdots)}\right)$. Then there exists a unique rational
function $q(v)$ such that $F$ is semiformally conjugate to
$$G_F=\left(\matrix{x+x^{k+1}p(v)+x^{2k+1}q(v)\cr v}\right)$$ in
$\overline \C\setminus \{p(v)=0\}$. In addition, $q(v)=\frac
{s(v)}{p(v)^{2k+1}}$ where $s(v)$ is a polynomial of degree
$2k+2+2k\partial(p(v))$.

\end{teo}

\proof Follows the same procedure to the proof of theorem
\ref{semiconj}. Define by induction the sequences  $F_n$ and $H_n$
from the initial condition  $F_0=F$ and $H_0=id$, and for $j>0$
($j\ne k$), if $$
F_{j-1}-G_F=\left(\matrix{x^{k+j+1}\phi_{1,j}(v)+\cdots\cr
x^{k+j}\phi_{2,j}(v)+\cdots}\right)$$
 define
$$H_j=(x,v)=\left(\matrix{x+x^{j+1}h_{1,j}(v)\cr v+x^j
h_{2,j}(v)}\right)\hbox{ and } F_j=H_j^{-1}\circ F_{j-1}\circ
H_j,$$ where $(h_{1,j},h_{2,j})$ is the unique solution of
(\ref{eq:dicrit}), i.e. $h_{1,j}=\frac
{\phi_{1,j}(v)}{(k-j)p(v)}+\frac {p'(v)\phi_{2,j}(v)}{jp(v)^2}$
and $h_{2,j}(v)=-\frac {\phi_{2,j}(v)}{jp(v)}$. See that the
degree of $h_{1,j}$ and $h_{2,j}$ are $j+1$ and $j+2$.

In the case $j=k$ we have that $$
F_{k-1}-G_F=\left(\matrix{x^{2k+1}(\phi_{1,k}(v)-q(v))+\cdots\cr
x^{2k}\phi_{2,k}(v)+\cdots}\right),$$ where $q(v)$ is the unique
function such that the system of equations  (\ref{eq:dicrit}) has
solution, i.e. $h_{2,k}(v)=-\frac {\phi_{2,k}(v)}{kp(v)}=\frac
{\phi_{1,k}(v)-q(v)}{p'(v)}$, therefore $$q(v)=\phi_{1,k}(v)-\frac
{p'(v)\phi_{2,k}(v)}{kp(v)}.$$ Finally,   it is easy to prove
inductively that $ \phi_{i,j}$ ($i=1,2$)  are of the form $\frac
{\psi_{i,j}}{p(v)^{2j}}$ where $ \psi_{1,j}$ and $ \psi_{2,j}$ are
polynomials of degree $k+j+1+2j\partial(p(v))$ and
$k+j+2+2j\partial(p(v))$ respectively, and then we conclude that
$q(v)$ is of form $q(v)=\frac {s(v)} {p(v)^{2k+1}}$ where $s(v)$
is a polynomial of degree $2k+2+2k\partial(p(v))$.\fin

{\bf Acknowledgmet:} The author are grateful to Prof. C\'esar
Camacho for suggesting the problem. This work was supported by the
IMPA-CNPq, Rio de Janeiro, Brasil.

\vskip 1cm \small {\bf Fabio Enrique Brochero Mart\'{\i}nez}\\
Instituto de Matem\'atica Pura e Aplicada (IMPA)
\\ Estrada Dona Castorina 110\\CEP 22460-320\\ Rio de Janeiro, RJ\\ {\bf Brasil}
\\ fbrocher@impa.br

\end{document}